\documentclass[11pt,twoside]{article}
\usepackage{amsthm,amsfonts,mathtools,amscd,amssymb, mathrsfs}
\usepackage{euscript}
\usepackage{enumitem}
\usepackage{graphicx}
\usepackage{tikz}
\usepackage[font=small]{caption}
\usepackage{subcaption}

\usepackage{cite}
\usepackage[utf8]{inputenc}
\usepackage[english]{babel}
\usepackage{xcolor}

\usepackage{jmpag-n} 

\hypersetup{linkcolor=darkblue, urlcolor=darkred, citecolor=darkgreen}

\begin{document}



\title[Approximate Controllability Problems for the Heat Equation in a Half-Plane]%
{Approximate Controllability Problems for the Heat Equation in a Half-Plane Controlled by the Dirichlet Boundary Condition with a Bounded Control}

\author{Larissa  Fardigola}
\address{B. Verkin Institute for Low Temperature Physics and Engineering of the
	National Academy of Sciences of Ukraine, 47 Nauky Ave., Kharkiv, 61103, Ukraine,\\
	V.N. Karazin Kharkiv National University, 4 Svobody Sq., Kharkiv, 61022, Ukraine}
\email{fardigola@ilt.kharkov.ua}

\author{Kateryna  Khalina}
\address{B. Verkin Institute for Low Temperature Physics and Engineering of the
	National Academy of Sciences of Ukraine, 47 Nauky Ave., Kharkiv, 61103, Ukraine}
\email{khalina@ilt.kharkov.ua}


\BeginPaper 


\newcommand{\LL}{L_{-1,(3)}^2(\mathbb R_+)}
\newcommand{\LLs}[1]{L_{#1,(3)}^2(\mathbb R_+)}
\newcommand{\fnl}{\left[\kern-0.1em\left]}
\newcommand{\fnr}{\right[\kern-0.1em\right]}
\newcommand{\sgn}{\mathop{\mathrm{sgn}}}
\newcommand{\supp}{\mathop{\mathrm{supp}}}
\newcommand{\zco}{{\text{$\bigcirc$\kern-6.4pt\raisebox{-0.5pt}{$0$}}}}
\newcommand{\R}{\mathbb R}
\newcommand{\N}{\mathbb N}


\setlength{\textfloatsep}{4pt plus 0.5pt minus 1.0pt}

\setlength{\floatsep}{4pt plus 0.5pt minus 1.0pt}

\setlength{\intextsep}{4pt plus 0.5pt minus 1.0pt}



\begin{abstract}
In the paper, the problems of approximate
controllability are  
studied for the control system $w_t=\Delta w$, $w(0,x_2,t)=u(x_2,t)$,
$x_1\in\R_+=(0,+\infty)$, $x_2\in\mathbb R$, $t\in(0,T)$, where $u$  is a control belonging to a special subset of $L^\infty(\R\times (0,T))\cap L^2(\R\times (0,T))$. It is proved that each initial state belonging to $L^2(\R_+\times\R)$ is approximately controllable to an arbitrary end state belonging to $L^2(\R_+\times\R)$ by applying these controls. A numerical algorithm of solving the approximate controllability problem for this system is given.
The results are illustrated by an example.

\key{heat equation, controllability, approximate controllability, half-plane}

\end{abstract}



\section{Introduction}\label{sect1}

We consider the heat equation  in a half-plane
\begin{align}
&w_t= \Delta w,&& x_1\in\R_+,\ x_2\in\mathbb{R},\ t\in(0,T),\label{eq}
\intertext{controlled by the Dirichlet boundary condition}
&w\big(0,(\cdot)_{[2]},t\big)=u((\cdot)_{[2]},t),&& x_2\in\mathbb{R},\ t\in(0,T),\label{bc}
\intertext{under the initial condition}
&w\big((\cdot)_{[1]},(\cdot)_{[2]},0\big)=w^0,&&x_1\in\R_+,\ x_2\in\mathbb{R}, \label{ic}
\end{align}
where $\R_+=(0,+\infty)$, $T>0$,  $\Delta =(\partial/\partial x_1)^2+(\partial/\partial x_2)^2$,
\begin{equation}
\label{contr}
u\in U[0,T]  
\end{equation}
is a  control, 
\begin{equation*}
U[0,T]=\left\{\varphi\in L^\infty(\R\times (0,T))\, \middle|\, \sup_{t\in[0,T]}|\varphi(\cdot,t)|\in L^2(\R) \right\}
\end{equation*}
is the set of admissible controls.
The subscripts $[1]$ and $[2]$ associate with the variable numbers, e.g. $(\cdot)_{[1]}$ and $(\cdot)_{[2]}$ correspond to $x_1$ and $x_2$, respectively, if we consider $f(x)$, $x\in\mathbb R^2$. This problem is considered in spaces of the Sobolev type (see details in Section \ref{sect2}).

Controllability problems for the heat equation were studied both in  bounded and unbounded domains. However, most of the papers studying these problems deal with domains bounded with respect to the spatial variables (see \cite{A, Em1, Em2, FatRus, FatRus1, Fat, AFOYI, WL, IY, LR, LTZua, MRT, Rus} and the references therein). 
Note that each initial state of the heat equation in a bounded domain which belongs to a Sobolev space of non-positive order can be driven to the origin in an arbitrarily small time by an $L^2$ boundary control (see, e.g. \cite{AFOYI,LR}).
Although there are quite a few papers considering domains unbounded with respect to the spatial variables  \cite{Barb, CabMenZua, CMV, LDCK, DWZh, FKh, FKh2, FKh3, FKh4, FKh5, GonTer, MenCab, MZua1, MZua2, Mill, DD, Ter, TerZua}, there are much less papers (see, e.g. \cite{FKh3,FKh6,MZua2}) where the boundary controllability of the heat equation was studied in a half-plain. 


For a bounded domain $\Omega\subset \R^n$ with the boundary $\partial \Omega$ of class $C^2$ (which is considered instead of the domain $\R_+\times\R$), it is well-known that the control system of the form \eqref{eq}--\eqref{ic} is null-controllable for a given time $T>0$. This result was obtained by using Carleman inequalities (see, e.g. \cite{AFOYI,LR}). 

For unbounded domains, the situation is essentially different. There exist pairs of initial and target states where the initial state can be driven to the end state by means of control system \eqref{eq}--\eqref{ic}, and there exist those where the initial state cannot be driven to the end state by means of this system. For instance, there is no initial data in any negative Sobolev space that may be driven to zero in finite time (see \cite{MZua2}). 

Although the statement of approximate null-controllability for control system \eqref{eq}--\eqref{ic} is considered known, there is no papers containing  direct and clear proof of it. E.g. in \cite{MZua2}, the authors  assert that ``one can easily prove the approximate controllability directly both in the case of bounded and unbounded domains'', but they  give neither a proof of this fact nor a reference for it.

In this paper we prove the approximate controllability for control system \eqref{eq}--\eqref{ic}. In particular, we construct controls solving the approximate controllability problem and give a numerical algorithm of solving this problem.


In \cite{MZua2}, the lack of the null-controllability property for the linear heat equation on the half-space $\R_+\times\R^n$ with
an $L^2$ Dirichlet boundary control is studied.  By rewriting system on the similarity variables and separating of variables, the multi-dimensional control problem is reduced to an infinite family of one-dimensional control systems in weighted Sobolev spaces. Due to results obtained for one-dimensional control system (see \cite{MZua1}), it was shown that controllable data have Fourier coefficients that grow exponentially for large frequencies. So, any initial state from Sobolev spaces cannot be driven to the origin by an $L^2$ Dirichlet boundary control in a time $T>0$.
In \cite{FKh3}, control system \eqref{eq}--\eqref{ic} was studied with a control of the form $u((\cdot)_{[2]},t)=\delta((\cdot)_{[2]})v(t)$ ($\delta$ is the Dirac distribution, $v\in L^\infty(0,T)$), which is not bounded in contrast to the present paper. In fact, the control $\delta((\cdot)_{[2]})v(t)$ is a distribution of the class $H^{-1}(\R)$ with respect to $x_2$. That is why the control system was considered in the Sobolev space $H^{-1}(\R_+\times\R)$ in \cite{FKh3} in contrast to the  space $H^0(\R_+\times\R)=L^2(\R_+\times\R)$ in the present paper. In addition, the spaces $\mathcal H^s$ and $\mathcal H_s$ based on functions of the form $f(|x|)$ were constructed and studied to investigate the controllability problems in \cite{FKh3} because the states of the control system which could be approximately steered to the origin had the form $w^0(x)=\frac\partial{\partial x_1}\omega^0(|x|)$. So, under the control of the form $u((\cdot)_{[2]},t)=\delta((\cdot)_{[2]})v(t)$,  the set of states which can be approximately driven to the origin is not coincide with $H^{-1}(\R_+\times\R)$ in contrast to the case of bounded controls in the present paper, where this set is equal to $H^0(\R_+\times\R)=L^2(\R_+\times\R)$.
In \cite{FKh3}, both necessary
and sufficient conditions for controllability and sufficient conditions for approximate controllability in a given time $T$ under a control $u$ bounded by a
given constant were obtained in terms of solvability of a Markov power moment problem. In addition, orthogonal bases were constructed in special spaces $\mathcal H^s$ and $\mathcal H_s$ of the Sobolev type in this paper. Using these bases, necessary and sufficient conditions for approximate controllability and numerical solutions to the approximate controllability problem were obtained in \cite{FKh3}. In \cite{FKh6}, the results of \cite{FKh3} were extended to the case of the Neumann boundary control, i.e. to the case where condition \eqref{bc} is replaced by the condition $w_{x_1}\big(0,(\cdot)_{[2]},t\big)=\delta((\cdot)_{[2]})v(t)$.
The boundary controllability of the wave equation in a half-plane $\R_+\times\mathbb R$ with a pointwise control on the bound  was studied in \cite{FLVJMPAG05, FLVJMPAG15, LVF0}.


In the present paper,  the approximate controllability problem for system \eqref{eq}--\eqref{ic} is studied in Sobolev spaces under controls from $U[0,T]$, in particular, $w^0\in L^2(\R_+\times\R)$ and $w(\cdot,t)\in L^2(\R_+\times\R)$, $t\in[0,T]$.
We show that $L^2(\R\times (0,T))$-controls are not appropriate to consider approximate controllability property for $w^0\in L^2(\R_+\times\R)$ and $w(\cdot,t)\in L^2(\R_+\times\R)$, $t\in[0,T]$, because there exists a control $u\in L^2(\R\times (0,T))$ (with compact supports) for which the end state $w(\cdot,T)$ of the solution to \eqref{eq}--\eqref{ic} does not belong to $L^2(\R_+\times\R)$ for the initial state $w^0=0\in L^2(\R_+\times\R)$ (see Example \ref{contrex} in Section \ref{sect3} below). That is why we consider the narrower set of controls $U[0,T]$, which provides the condition $w(\cdot,t)\in L^2(\R_+\times\R)$, $t\in[0,T]$, for any $w^0\in L^2(\R_+\times\R)$ (see Theorem \ref{th-sol} in Section \ref{sect3} below). Roughly speaking, we consider a specific subset of bounded controls in $L^2(\R\times (0,T))$. We prove that each initial state $w^0\in L^2(\R_+\times\R)$ of system \eqref{eq}--\eqref{ic} can be driven to an arbitrary neighbourhood of any target state $w^T\in L^2(\R_+\times\R)$ by choosing an appropriate control $u\in U[0,T]$, in other words, a state $w^0\in L^2(\R_+\times\R)$ is approximately controllable to a target state $w^T\in L^2(\R_+\times\R)$ in a given time $T$ (see Theorem \ref{thappc} in Section \ref{sect4} below). The method of proving this assertion is constructive. This allows to provide a numerical algorithm of solving the approximate controllability problem for system \eqref{eq}--\eqref{ic}. To this aid, we consider the odd extension of $w$ and $w^0$ with respect to $x_1$ 
and obtain control problem \eqref{eq1}, \eqref{ic1} in Section \ref{sect3}. Then we develop the state and the control in this new system in the Fourier series with respect to a basis generated by Hermite functions that allows us to reduce the 2-d problem to a family of the 1-d ones. Since we consider the approximate controllability, we can solve only a finite number of these 1-d problems. To construct controls solving them, we apply the method introduced in \cite{FKh} for solving the approximate controllability problem for the 1-d heat equation controlled by the Dirichlet boundary condition.
We should note that similar development was used in \cite{MZua2} to reduce the 2-d problem to a family of the 1-d ones, but the basis of the eigenfunctions of the differential operator had been obtained after using similarity variables and weighted Sobolev spaces. In our paper, we use developing in the Fourier series directly in $L^2(\R^2)$ without using similarity variables that simplifies the numerical method. Moreover, we apply the Fourier transform and its inverse to analyze the solution to control problem.

The paper is organized as follows:
in Section \ref{sect2}, some notations and definitions are given; in Section \ref{sect3},  the controllability problem is formulated for the control system \eqref{eq}--\eqref{ic}, and  preliminary results are given;
in Section \ref{sect4}, the approximate controllability results are obtained;
in Section \ref{sect5}, a numerical algorithm of solving the approximate controllability problem for system \eqref{eq}--\eqref{ic} is given; in Section \ref{ex}, the results are illustrated by an example.


\section{Notation}\label{sect2}

Let us introduce the spaces used in the paper. Let $n\in\mathbb N$.
By $|\cdot|$, we denote the Euclidean norm in $\mathbb{R}^n$.

Let $\mathscr{S}(\mathbb R^n)$ be the Schwartz space of rapidly decreasing functions \cite{Schw}, $\mathscr{S}'(\mathbb{R}^n)$ be the dual space. 

Let $D=\big(-i\partial/\partial x_1,\ldots,-i\partial/\partial x_n\big)$, $D^\alpha=\big(-i(\partial/\partial x_1)^{\alpha_1},\ldots,-i(\partial/\partial x_n)^{\alpha_n}\big)$, where $\alpha=(\alpha_1,\ldots,\alpha_n)\in\N_0^n$ is multi-index, $\N_0=\N\cup\{0\}$.

For $s=\overline{0,2}$, consider
\begin{equation*}
H^s(\mathbb{R}^n)=\left\{\varphi\in L^2(\mathbb{R}^n)\,\middle|\, \forall\alpha\in\mathbb{N}_0^n \quad \big(\alpha_1+\cdots+\alpha_n\leq s \Rightarrow D^\alpha\varphi\in L^2(\mathbb{R}^n)\big)\right\}
\end{equation*}
with the norm
\begin{equation*}
\left\| \varphi\right\|^s=\left(\sum_{\alpha_1+\cdots+\alpha_n\leq s} \frac{ s!}{(s-|\alpha|)!\alpha!} \left(\left\| D^\alpha\varphi\right\|_{L^2(\mathbb{R}^n)} \right)^2\right)^{1/2},\quad\varphi\in H^s(\mathbb{R}^n),
\end{equation*}
and $H^{-s}(\mathbb{R}^n)=\left( H^s\right(\mathbb{R}^n))^*$ with the  norm $\left\|\cdot\right\|^{-s}$ associated with the strong topology of the adjoint space. We also have $H^0(\mathbb{R}^n)=L^2(\mathbb{R}^n)=\left( H^0\left(\mathbb{R}^n\right)\right)^*$.

For $m=\overline{-2,2}$, consider
\begin{equation*}
H_m(\mathbb{R}^n)=\left\{ \psi\in L_{\text{loc}}^2(\mathbb{R}^n)\,\middle|\, \left(1+|\sigma|^2 \right)^{m/2}\psi\in L^2(\mathbb{R}^n)\right\}
\end{equation*}
with the norm
\begin{equation*}
\left\| \psi\right\|_m=\left\|\left(1+|\sigma|^2 \right)^{m/2}\psi \right\|_{L^2(\mathbb{R}^n)},\quad\psi\in H_m(\mathbb{R}^n).
\end{equation*}
Evidently, $H_{-m}(\mathbb{R}^n)=\left( H_m(\mathbb{R}^n)\right)^*$.

Let $\langle f,\varphi\rangle$ be the value of a distribution $f\in \mathscr S'(\mathbb R^n)$ on a test function $\varphi \in \mathscr S(\mathbb R^n)$. 

By $\mathscr{F}: \mathscr S'(\mathbb R^n)\to \mathscr S'(\mathbb{R}^n)$ denote the Fourier transform operator
with the domain $\mathscr S'(\mathbb R^n)$. This operator is an extension of the
classical Fourier transform operator which is an isometric
isomorphism of $L^2(\mathbb{R}^n)$. The extension is given by the formula
\begin{equation*}
\langle \mathscr{F} f,\varphi\rangle=\langle f,\mathscr{F}^{-1}\varphi\rangle,\quad f\in \mathscr S'(\mathbb R^n),\ \varphi\in \mathscr S(\mathbb R^n).
\end{equation*}
The operator $\mathscr{F}$ is an isometric isomorphism of $H^m(\mathbb{R}^n)$ and $H_m(\mathbb{R}^n)$, $m=\overline{-2,2}$, \cite[Chap.~1]{VG}.

A distribution $f\in \mathscr S'\left(\mathbb R^2\right)$  is said  to be
\emph{odd with respect to} $x_1$, if for all $\varphi\in \mathscr S(\mathbb R^2)$, we have
$\big\langle f,\varphi\big((\cdot)_{[1]},(\cdot)_{[2]}\big)\big\rangle=-\big\langle f,\varphi\big(-(\cdot)_{[1]},(\cdot)_{[2]}\big)\big\rangle$.

Let $m=\overline{-2,2}$. By $\widetilde{H}^m\left(\mathbb{R}^2\right)$ \Big(or $\widetilde{H}_m\left(\mathbb{R}^2\right)$\Big), denote the subspace of all distributions in
$H^m\left(\mathbb{R}^2\right)$ \Big(or $H_m\left(\mathbb{R}^2\right)$, respectively\Big) that are odd with respect to $x_1$. Evidently, $\widetilde{H}^m\left(\mathbb{R}^2\right)$ \Big(or $\widetilde{H}_m\left(\mathbb{R}^2\right)$\Big) is a closed
subspace of $H^m\left(\mathbb{R}^2\right)$ \Big(or $H_m\left(\mathbb{R}^2\right)$, respectively\Big).

For $s=\overline{0,2}$, consider
\begin{align*}
H_\zco^s=&\left\{\varphi\in L^2(\mathbb{R}_+\times\mathbb{R})\,\middle|\, \Big( \forall\alpha\in\mathbb{N}_0^2\quad \big(\alpha_1+\alpha_2\leq s\Rightarrow D^\alpha\varphi\in L^2(\mathbb{R}_+\times\mathbb{R})\big)\Big)
	\right.\\
	&\kern30ex\left.
\wedge  \Big(\forall k=\overline{0,s-1} \quad  D^{(k,0)}\varphi(0^+,(\cdot)_{[2]})=0\Big)\right\}
\end{align*}
with the norm
\begin{equation*}
\left\| \varphi\right\|_\zco^s=\left(\sum_{\alpha_1+\alpha_2\leq s} \frac{ s!}{(s-|\alpha|)!\alpha!} \left(\left\| D^\alpha\varphi\right\|_{L^2(\mathbb{R}_+\times\mathbb{R})} \right)^2\right)^{1/2},\quad \varphi\in H_\zco^s\,,
\end{equation*} 
and $H_\zco^{-s}=\left( H_\zco^s\right)^*$ with the  norm 
$\left\|\cdot\right\|_\zco^{-s}$ associated with the strong topology of the adjoint space.
We have
$$
H_\zco^0=L^2(\mathbb{R}_+\times\mathbb{R}).
$$

\begin{remark} 
	Let $\varphi \in H_\zco^s$\,, $s=\overline{0,2}$. Let $\widetilde\varphi$ be 
	its odd extension with respect to $x_1$, i.e.
	$\widetilde\varphi(x_1,x_2)=\varphi(x_1,x_2)$ if $x_1\ge0$ and  $\widetilde\varphi(x_1,x_2)=-\varphi(-x_1,x_2)$ otherwise. Then $\widetilde\varphi \in \widetilde H^s\left(\mathbb{R}^2\right)$, $s=\overline{0,2}$.
	The converse assertion is true for $s=0,1$, and it is not true for $s=2$. That is why the odd extension with respect to $x_1$ of a distribution $f\in H_\zco^{-2}$
	may not belong to $\widetilde{H}^{-2}\left(\mathbb{R}^2\right)$. However, the following theorem holds.
\end{remark}

\begin{theorem}
	\label{todext}
	Let $f\in H_\zco^0$ and there exists $f\big(0^+,(\cdot)_{[2]}\big)\in H^0(\mathbb{R})$. Then $f_{x_1x_1}\in
	H_\zco^{-2}$ can be extended to a distribution $F\in \widetilde{H}^{-2}\left(\mathbb{R}^2\right)$ such that $F$ is odd with respect to $x_1$. This distribution is given by the formula
	\begin{equation}
	\label{oddext}
	F=\widetilde f_{x_1x_1}-2f(0^+,\big(\cdot)_{[2]}\big)\delta',
	\end{equation}
	where $\widetilde f$ is the odd extension of $f$ with respect to $x_1$, $\delta$ is the Dirac distribution with respect to $x_1$.
\end{theorem}

In the case $f\in H_\zco^{-1/2}$, corresponding theorem has been proved in \cite{LVF0}. The proof
of Theorem \ref{todext}  is analogous to the proof of the mentioned theorem.

\section{Problem formulation and preliminary results}\label{sect3}

We consider control system \eqref{eq}--\eqref{ic} in $H_\zco^{-l}$, $l=\overline{0,2}$, i.e. $\left(\frac
d{dt}\right)^s w:[0,T]\to H_\zco^{-2s}$, $s=0,1$, $w^0\in H_\zco^0$\,. 

Let $w^0,w(\cdot,t)\in H_\zco^0$\,, $t\in [0,T]$. Let $W^0$ and $W(\cdot,t)$ be the odd extensions of $w^0$ and $w(\cdot,t)$ with
respect to $x_1$, respectively, $t\in [0,T]$. If $w$ is a solution to control system \eqref{eq}--\eqref{ic}, then $W$ is a solution to control system
\begin{align}
&W_t=\Delta W-2u((\cdot)_{[2]},t)\delta',\quad  t\in(0,T),
\label{eq1}
\\
&W\big((\cdot)_{[1]},(\cdot)_{[2]},0\big)=W^0 
\label{ic1}
\end{align}
according to Theorem  \ref{todext}.  Here $\left(\frac
d{dt}\right)^s W:[0,T]\to \widetilde{H}^{-2s}\left(\mathbb{R}^2\right)$, $s=0,1$, $W^0\in \widetilde{H}^0\left(\mathbb{R}^2\right)$. The converse assertion
is also true. Let $W^0, W(\cdot,t)\in \widetilde{H}^0\left(\mathbb{R}^2\right)$, $t\in [0,T]$. Let $w^0$ and $w(\cdot,t)$ be the restrictions of $W^0$ and $W(\cdot,t)$ to $(0,+\infty)$ with
respect to $x_1$, respectively, $t\in [0,T]$. If $W$ is a solution to \eqref{eq1}, \eqref{ic1},
then $w$ is a solution to \eqref{eq}--\eqref{ic} because 
\begin{equation}
\label{bc1}
W\big(0^+,(\cdot)_{[2]},t\big)=u((\cdot)_{[2]},t)\quad\text{for almost all}\ t\in[0,T]
\end{equation}
according to Lemma \ref{lempoint} (see below). Assume that $w^T\in H_\zco^0$\,. Evidently, $w\big((\cdot)_{[1]},(\cdot)_{[2]},T\big)=w^T$  iff  $W\big((\cdot)_{[1]},(\cdot)_{[2]},T\big)=W^T$.
Here $W^T$ is the odd extension of $w^T$ with respect to $x_1$ and $W^T\in \widetilde{H}^0\left(\mathbb{R}^2\right)$.

Thus, control systems \eqref{eq}--\eqref{ic} and \eqref{eq1}, \eqref{ic1} are equivalent. Therefore, basing on this reason, we will further consider control system  \eqref{eq1}, \eqref{ic1} instead of original system \eqref{eq}--\eqref{ic}.

Let $T>0$, $W^0\in\widetilde{H}^0\left(\mathbb{R}^2\right)$.
By $\mathscr{R}_T\left(W^0\right)$, denote the set of all states $W^T\in\widetilde{H}^0\left(\mathbb{R}^2\right)$ for which there exists a control $u\in U[0,T]$ such that there exists a unique solution $W$ to system \eqref{eq1}, \eqref{ic1} such that  $W\big((\cdot)_{[1]},(\cdot)_{[2]},T\big)=W^T$. 

\begin{definition}
	A state $W^0\in\widetilde{H}^0\left(\mathbb{R}^2\right)$ is said to be controllable to a target state $W^T\in\widetilde{H}^0\left(\mathbb{R}^2\right)$ in a given time $T>0$ if
	$W^T\in\mathscr{R}_T\left(W^0\right)$.
\end{definition}

In other words, a state
$W^0\in\widetilde{H}^0\left(\mathbb{R}^2\right)$ is  controllable to a target state $W^T\in\widetilde{H}^0\left(\mathbb{R}^2\right)$ in a given time $T>0$ if there exists a control $u\in U[0,T]$ such that there exists a unique solution $W$ to system \eqref{eq1}, \eqref{ic1} and $W\big((\cdot)_{[1]},(\cdot)_{[2]},T\big)=W^T$.

\begin{definition}
	\label{def-appr}
	A state $W^0\in\widetilde{H}^0\left(\mathbb{R}^2\right)$ is said to be approximately controllable to a target state $W^T\in\widetilde{H}^0\left(\mathbb{R}^2\right)$ in a given time $T>0$ if
	$W^T\in\overline{\mathscr{R}_T\left(W^0\right)}$, where the closure is considered in the space $\widetilde{H}^0\left(\mathbb{R}^2\right)$.
\end{definition}

In other words, a state $W^0\in\widetilde{H}^0\left(\mathbb{R}^2\right)$ is approximately
controllable  to a target state $W^T\in\widetilde{H}^0\left(\mathbb{R}^2\right)$ in a given time
$T>0$ if for each $\varepsilon>0$, there exists $u_\varepsilon\in  U[0,T]$ such that there exists a unique solution $W_\varepsilon$ to
system \eqref{eq1}, \eqref{ic1} with $u=u_\varepsilon$ and $\left\|
W_\varepsilon\big((\cdot)_{[1]},(\cdot)_{[2]},T\big)-W^T \right\|^0<\varepsilon$.

Using the fundamental solution to the heat operator (see, e.g. \cite[Chapter 7]{SHB}),  we obtain  the unique solution to system \eqref{eq1},  \eqref{ic1} 
\begin{equation}
\label{sol1}
W(x,t)=\mathcal W_0(x,t)+\mathcal W_u(x,t), \quad x\in\mathbb{R}^2,\ t\in[0,T],
\end{equation}
where
\begin{align}
\label{wo}
\mathcal W_0(x,t)&=\frac{1}{4\pi t}e^{-|x|^2/(4t)}*W^0(x),\quad x\in\mathbb{R}^2,\ t\in[0,T],
\\
\label{wu}
\mathcal W_u(x,t)&=\frac{x_1}{\pi}\int_0^t \frac{1}{4\xi^2}e^{-|x|^2/(4\xi)}*_{[2]} u(x_2,t-\xi)\,d\xi,\quad x\in\mathbb{R}^2,\ t\in[0,T].
\end{align}
Here $*$ is the convolution with respect to $x$ and $*_{[2]}$ is the convolution with respect to $x_2$.
\begin{theorem}
	\label{th-sol}
	Let $u\in U[0,T]$, 
	$W^0\in\widetilde{H}^0\left(\mathbb{R}^2\right)$. Then,
	\begin{enumerate}
		\renewcommand{\labelenumi}{\parbox{0.7cm}{\upshape (\roman{enumi})\hfill}}
		\renewcommand{\theenumi}{(\roman{enumi})}
		\item \label{th-sol-i}
		$\mathcal W_0(\cdot,t)\in\widetilde H^0\left(\mathbb R^2\right)$, $t\in[0,T];$
		\item \label{th-sol-ii}
		$\mathcal W_0(\cdot,t)\in C^\infty\left(\mathbb R^2\right)$, $t\in(0,T];$
		\item \label{th-sol-iii}
		$\|\mathcal W_u(\cdot,t)\|^0\leq 2t^{1/4} 
		\left\|\sup_{t\in[0,T]}|u(\cdot,t)|\right\|_{L^2(\R)}$, $t\in(0,T].$
	\end{enumerate}
\end{theorem}

\begin{proof}
	Denote $V^0=\mathscr{F} W^0$, 
	$\mathcal V_0(\cdot,t)=\mathscr{F}_{x\to\sigma} \mathcal W_0(\cdot,t)$, 
	We have
	\begin{equation}
	\label{a2}
	\mathcal V_0(\sigma,t) = e^{-t|\sigma|^2} V^0(\sigma),\quad\sigma\in\mathbb R^2,\ t\in[0,T].
	\end{equation}
	Therefore,
	\begin{equation}
	\label{a4}
	\|\mathcal W_0(\cdot,t)\|^0=\|\mathcal V_0(\cdot,t)\|_0
	\leq \|V^0\|_0=\|W^0\|^0,\quad t\in[0,T],
	\end{equation}
	i.e. \ref{th-sol-i}  holds.
	
	Let $\alpha=(\alpha_1,\alpha_2)\in\N_0^2$. We have
	\begin{align*}
	\left|\sigma_1^{\alpha_1}\sigma_2^{\alpha_2} \mathcal V_0(\sigma,t) \right|^2\leq (1+|\sigma|^2)^{1+\alpha_1+\alpha_2}e^{-2t|\sigma|^2}(1+|\sigma|^2)^{-1}&|V^0(\sigma)|^2,
	\\ 
	&\sigma \in \mathbb R^2,\ t\in(0,T].
	\end{align*}
	Since 
	\begin{equation*}
	\xi^me^{-\beta\xi}\leq\left(\frac{m}{\beta e}\right)^m,\quad \xi\geq 0,\ \beta>0,\ m\in\mathbb N_0,
	\end{equation*}
	then
	\begin{align}
	\label{a5}
	\left\| D^\alpha \mathcal W_0(\cdot,t) \right\|^0
	&=\left\|(\cdot)_{[1]}^{\alpha_1}(\cdot)_{[2]}^{\alpha_2} \mathcal V_0(\cdot,t) \right\|_0
	\leq e^t\left(\frac{1+\alpha_1+\alpha_2}{2te}\right)^{(1+\alpha_1+\alpha_2)/2}\big\|V^0\big\|_{-1}
	\notag\\
	&\leq e^t\left(\frac{1+\alpha_1+\alpha_2}{2te} \right)^{(1+\alpha_1+\alpha_2)/2}\|V^0\|_0
	\notag\\
	&=e^t\left(\frac{1+\alpha_1+\alpha_2}{2te}\right)^{(1+\alpha_1+\alpha_2)/2}\big\|W^0\big\|^0,\quad  t\in(0,T],
	\end{align}
	i.e.  \ref{th-sol-ii}  holds.
	
	Put
	$$
	g(x_2)= \sup_{t\in[0,T]}|u(x_2,t)|,\quad x_2\in\R.
	$$
	Since $u\in U[0,T]$, we have $g\in L^2(\R)$. Put
	$G=\mathscr F g$. Then we have
	\begin{equation}
	\label{ab1}
	\left\|\sup_{t\in[0,T]}|u(\cdot,t)|\right\|_{L^2(\R)} 
	=\|g\|_{L^2(\R)}=\|G\|_{L^2(\R)}.
	\end{equation}
	Taking into account \eqref{wu}, we get
	$$
	\mathcal W_u(x,t)
	=\frac{x_1}{\pi}\int_0^t \frac{e^{-x_1^2/(4\xi)}}{4\xi^2} 
	\int_{-\infty}^\infty e^{-(x_2-\mu)^2/(4\xi)} u(\mu,t-\xi)\, d\mu\, d\xi,\  x\in\R^2,\ t\in[0,T].
	$$ 
	Therefore,
	\begin{align}
	\label{ab2}
	&|\mathcal W_u(x,t)|
	\leq 
	\frac{|x_1|}{\pi}\int_0^t \frac{e^{-x_1^2/(4\xi)}}{4\xi^2} 
	\int_{-\infty}^\infty e^{-(x_2-\mu)^2/(4\xi)} g(\mu)\, d\mu\, d\xi
	\notag\\
	&\ \ 
	=\frac{|x_1|}{\pi}\int_0^t \frac{1}{4\xi^2}e^{-|x|^2/(4\xi)}*_{[2]} g(x_2)\,d\xi
	=|\mathcal W_g(x,t)|,\quad x\in\R^2,\ t\in[0,T].
	\end{align}
	Here, we have again applied \eqref{wu} to the last equality. Setting $\mathcal V_g(\cdot,t)= \mathscr F \mathcal W_g(\cdot,t)$, $t\in[0,T]$, we obtain
	$$
	\mathcal V_g(\sigma,t)
	=-\sqrt{\frac2\pi}i\sigma_1\int_0^t e^{-\xi|\sigma|^2} G(\sigma_2)\, d\xi,\quad \sigma \in\R^2,\ t\in[0,T].
	$$
	Then, for $t\in[0,T]$, we have
	\begin{align*}
	\left(\|\mathcal V_g(\cdot,t)\|_0\right)^2
	=\frac2\pi \int_0^t \int_0^t 
	\int_{-\infty}^\infty \sigma_1^2 e^{-(\xi+\mu)\sigma_1^2}\, d\sigma_1
	\int_{-\infty}^\infty  e^{-(\xi+\mu)\sigma_2^2} |G(\sigma_2)|^2\, d\sigma_2\, d\xi\, d\mu.
	\end{align*}
	Since
	$$
	\int_{-\infty}^\infty y^2 e^{-\alpha y^2}\,dy
	=-\frac{d}{d\alpha}
	\int_{-\infty}^\infty e^{-\alpha y^2}\,dy
	=-\left(\sqrt{\frac\pi\alpha}\right)'
	=\frac{\sqrt\pi}{2 \alpha^{3/2}}, \quad \alpha>0, 
	$$
	we obtain
	\begin{align*}
	\left(\|\mathcal V_g(\cdot,t)\|_0\right)^2
	&\leq \frac1{\sqrt\pi} \left(\|G\|_{L^2(\R)}\right)^2 
	\int_0^t \int_0^t \frac1{(\xi+\mu)^{3/2}}\,d\xi\,d\mu
	\\
	&= \frac2{\sqrt\pi} \left(\|G\|_{L^2(\R)}\right)^2 
	\int_0^t \left(\frac1{\sqrt\mu}-\frac1{\sqrt{\mu+t}}\right)d\mu
	\\
	&=\frac{4\left(2-\sqrt2\right)}{\sqrt\pi}\sqrt t \left(\|G\|_{L^2(\R)}\right)^2 ,\quad t\in[0,T].
	\end{align*}
	Taking into account \eqref{ab1} and \eqref{ab2}, we conclude that \ref{th-sol-iii}  holds.
\end{proof}

It follows from Theorem \ref{th-sol}\ref{th-sol-iii} that $\mathcal W_u(\cdot,T)\in \widetilde{H}^0\left(\R^2\right)\subset L^2\left(\R^2\right)$ if $U[0,T]\subset L^2\left(\R\times[0,T]\right) \cap L^\infty\left(\R\times[0,T]\right)$. The following example shows that the boundedness of a control $u$  plays a significant role for $\mathcal W_u(\cdot,T)$ to belong to $L^2\left(\R^2\right)$.

\begin{example}
	\label{contrex}
	Consider
	$$
	u(x_2,t)=\frac1{(T-t)^{5/8}}\Big(H\Big(x_2+2\sqrt{T-t}\Big) -H\Big(x_2-2\sqrt{T-t}\Big)\Big),
	\  x_2\in\R,\ t\in(0,T],
	$$ 
	where $H$ is the Heaviside function ($H(\nu)=0$ if $\nu<0$ and $H(\nu)=1$ otherwise).
	Evidently, $u\notin L^\infty\left(\R\times[0,T]\right)$, but $u\in L^2\left(\R\times[0,T]\right)$ because
	$$
	\left(\|u\|_{L^2\left(\R\times[0,T]\right)}\right)^2
	=\int_0^T \frac1{\xi^{5/4}}
	\int_{-2\sqrt\xi}^{2\sqrt\xi}dx_2\,d\xi
	=4\int_0^T \xi^{-3/4}\, d\xi=16T^{1/4}<\infty. 
	$$
	
	Let us show that $\mathcal W_u(\cdot,T)\notin  L^2\left(\R^2\right)$. With regard to \eqref{wu}, we get
	\begin{equation}
	\label{enc0}
	\mathcal W_u(x,T)
	=\frac{x_1}{\pi} \int_0^T \frac{e^{-|x|^2/(4\xi)}}{4\xi^{2+5/8}} *_{[2]}\Big(H\Big(x_2+2\sqrt\xi\Big) -H\Big(x_2-2\sqrt\xi\Big)\Big)\, d\xi,\quad x\in\R^2.
	\end{equation}
	First, we calculate and estimate the convolution in the right-hand side of \eqref{enc0}. We have
	\begin{align*}
	e^{-x_2^2/(4\xi)}  *_{[2]}\Big(H\Big(x_2+2\sqrt\xi\Big) &-H\Big(x_2-2\sqrt\xi\Big)\Big)
	=\int_{x_2-2\sqrt\xi}^{x_2+2\sqrt\xi} e^{-\mu^2/(4\xi)}\,d\mu
	\\
	&=2\sqrt\xi \int_{x_2/(2\sqrt\xi)-1}^{x_2/(2\sqrt\xi)+1} e^{-y^2}\, dy,\quad x_2\in\R, \xi\in(0,T],
	\end{align*}
	where $y=\mu/(2\sqrt\xi)$.
	Since
	$$
	\left|\frac{x_2}{2\sqrt\xi}\pm 1\right|^2 
	\leq \left(\frac{|x_2|}{2\sqrt\xi}+1\right)^2
	\leq 2\left(\frac{x_2^2}{4\xi}+1\right)
	=\frac{x_2^2}{2\xi}+2,\quad  x_2\in\R, \xi\in(0,T],
	$$
	we get
	\begin{align*}
	\left|e^{-x_2^2/(4\xi)} *_{[2]}\Big(H\Big(x_2+2\sqrt\xi\Big) -H\Big(x_2-2\sqrt\xi\Big)\Big)\right|
	&\geq \frac{4\sqrt\xi}{e^2} e^{-x_2^2/(2\xi)},
	\\ 
	&\qquad
	x_2\in\R, \xi\in(0,T].
	\end{align*}
	Taking this into account, we obtain from \eqref{enc0} that
	$$
	\big|\mathcal W_u(x,T)\big|
	\geq \frac{2|x_1|}{e^2\pi} 
	\int_0^T \frac{e^{-(x_1^2+2x_2^2)/(4\xi)}}{\xi^{1/8}} \frac{d\xi}{2\xi^2}
	\geq \frac{2|x_1|}{e^2\pi}
	\int_0^T \frac{e^{-|x|^2/(2\xi)}}{\xi^{1/8}} \frac{d\xi}{2\xi^2},\quad x\in\R^2.
	$$
	Replacing $|x|^2/(2\xi)$ by $y$ in the integral, we get
	\begin{equation}
	\label{enc1}
	\big|\mathcal W_u(x,T)\big|
	\geq \frac{2^{9/8}}{e^2\pi} \frac{|x_1|}{|x|^{9/4}}
	\int_{|x|^2/(2T)}^\infty e^{-y} y^{1/8}\, dy, \quad x\in\R^2.
	\end{equation}
	Setting
	$$
	F=\int_0^\infty e^{-y} y^{1/8}\, dy,
	$$
	we conclude that there exists $\varepsilon>0$ such that
	$$
	\int_{|x|^2/(2T)}^\infty e^{-y} y^{1/8}\, dy\geq \frac F{2^{9/8}}, \quad |x|\leq\varepsilon.
	$$
	Therefore \eqref{enc1} yields 
	$$
	\big|\mathcal W_u(x,T)\big|
	\geq \frac F{e^2\pi} \frac{|x_1|}{|x|^{9/4}}, \quad |x|\leq\varepsilon.
	$$
	Hence,
	$$
	\iint_{\R^2} \big|\mathcal W_u(x,T)\big|^2\, dx
	\geq  \frac {F^2}{e^4\pi^2} \iint_{|x|\leq\varepsilon}\frac{|x_1|^2}{|x|^{9/2}}\, dx=\infty,
	$$
	i.e. $\mathcal W_u(\cdot,T)\notin L^2\left(\R^2\right)$.
\end{example}

This example demonstrates that the set of admissible controls $U[0,T]$ being considered in the present paper cannot be extended to $L^2\left(\R\times[0,T]\right)$ if we want each control from the set of admissible controls to generate the end state belonging to $L^2\left(\R^2\right)$.

According to \eqref{sol1},  we have
\begin{align}
\label{re0}
\mathscr{R}_T(W^0)&=\left\{ W^T\in\widetilde{H}^0\left(\mathbb{R}^2\right) \,\middle|\,\exists u\in U[0,T]\quad
W^T=\mathcal W_0(\cdot,T)+\mathcal W_u(\cdot,T)\right\},
\end{align}
in particular,
\begin{equation}
\label{re1}
\mathscr{R}_T(0)=\left\{ W^T\in\widetilde{H}^0\left(\mathbb{R}^2\right) \,\middle|\,\exists u\in U[0,T]\quad W^T=\mathcal W_u(\cdot,T)\right\}.
\end{equation}
Taking into  account \eqref{sol1} and Theorem \ref{th-sol}, we obtain the following theorem.
\begin{theorem} 
	\label{reachprop}
	Let $T>0$. We have $f\in\mathscr{R}_T(W_0)$ iff $\left(f-\mathcal W_0(\cdot,T)\right)\in\mathscr{R}_T(0)$. We also have $f\in\overline{\mathscr{R}_T(W_0)}$ iff $\left(f-\mathcal W_0(\cdot,T)\right)\in\overline{\mathscr{R}_T(0)}$.
\end{theorem}

\begin{lemma}
	\label{lempoint}
	Let $u\in U[0,T]$ and $W^0\in \widetilde{H}^0\left(\mathbb{R}^2\right)$, $t\in [0,T]$. Let $W$ be a solution to \eqref{eq1}, \eqref{ic1}. Then \eqref{bc1} holds.
\end{lemma}

\begin{proof}
	According to Theorem \ref{th-sol}\ref{th-sol-ii}, $\mathcal W_0(\cdot,t)$ is  continuous on $\mathbb R^2$ for each $t\in(0,T]$. Moreover, $\mathcal W_0(\cdot,t)$ is odd with respect to $x_1$, $t\in[0,T]$. Hence,
	\begin{equation}
	\label{aa1}
	\mathcal W_0(0^+,(\cdot)_{[2]},t)=0,\quad t\in(0,T].
	\end{equation}
	
	Let us calculate $\mathcal W_u(0^+,(\cdot)_{[2]},t)$, $t\in[0,T]$. 
	It follows from \eqref{wu} that
	$$
	\mathcal W_u(x,t)=\frac{x_1}{\pi} \int_{-\infty}^{\infty} \int_0^t \frac1{4\xi^2} e^{-(x_1^2+\mu_2^2)/(4\xi)} u(x_2-\mu_2,t-\xi)\,d\xi\,d\mu_2,\  x\in\R^2,\ t\in[0,T].
	$$ 
	Changing variables in the integral:
	$\xi=\frac{x_1^2}{4\zeta^2}(1+\eta^2)$, $ \mu_2=|x_1|\eta$,
	we get 
	\begin{equation}
	\label{int1}
	\mathcal W_u(x,t)=\frac{\sgn x_1}{\pi} \int_{-\infty}^{\infty} \int_0^\infty
	F(x,t,\zeta,\eta)\, d\zeta\,d\eta,
	\end{equation}
	where
	\begin{align*}
	F(x,t,\zeta,\eta)=H\left(\!\zeta-\frac{|x_1|}{2\sqrt t}\sqrt{1+\eta^2}\right) \frac{2\zeta}{1+\eta^2} &e^{-\zeta^2} u \left(\! x_2-|x_1|\eta, t-\frac{x_1^2}{4\zeta^2}\left(1+\eta^2\right)\!\right), \\
	&
	x\in\R^2, \ t\in(0,T],\ \zeta\in\R_+,\ \eta\in\R.
	\end{align*} 
	With regard to \eqref{contr},  we obtain
	\begin{align}
	\label{est11}
	|F(x,t,\zeta,\eta)|\leq \frac{2\zeta}{1+\eta^2} e^{-\zeta^2}&\|u\|_{L^\infty(\R\times(0,T))} 
	\notag\\
	&\text{for a.a.}\ x\in\R^2, \ t\in(0,T],\ \zeta\in\R_+,\ \eta\in\R.
	\end{align}
	On the other hand, for a.a. $x_2\in\R$, $t\in(0,T]$, $\zeta\in\R_+$, $ \eta\in\R$, we have
	\begin{equation}
	\label{lim0}
	F(x,t,\zeta,\eta)\to \frac{2\zeta}{1+\eta^2} e^{-\zeta^2} u(x_2,t)\quad \text{as } x_1\to0^+.
	\end{equation}
	Taking into account \eqref{est11} and \eqref{lim0},  passing to the limit as $x_1\to0^+$, and applying Lebesgue's dominated convergence theorem to the right-hand side of \eqref{int1}, we obtain
	\begin{align}
	\label{aa3}
	\mathcal W_u\big(0^+,(\cdot)_{[2]},t\big)
	&=\frac1\pi u(x_2,t)\int_{-\infty}^{\infty} \frac{d\eta}{1+\eta^2} \int_0^\infty 2\zeta e^{-\zeta^2}\,d\zeta
	\notag\\
	&= u((\cdot)_{[2]},t),\quad\text{for a.a.}\ x_2\in\R,\ t\in[0,T].
	\end{align}
	Taking into account \eqref{sol1}, \eqref{aa1}, and \eqref{aa3}, we obtain  \eqref{bc1}.
\end{proof}


\section{Approximate controllability}
\label{sect4}

In this section we prove the main theorem.
\begin{theorem}
	\label{thappc}
	Each state $W^0\in\widetilde{H}^0\left(\mathbb{R}^2\right)$ is approximately controllable to a state $W^T\in\widetilde{H}^0\left(\mathbb{R}^2\right)$
	in a given time $T>0$, i.e. $\widetilde{H}^0\left(\mathbb{R}^2\right)=\overline{\mathscr{R}_T(W^0)}$.
\end{theorem}

First, consider the Hermite polynomials \cite[18.5.5 (and Tables 18.3.1, 18.5.1), 18.5.13]{NIST} 
\begin{equation}
\label{hermite}
\mathcal H_n(\mu)=(-1)^ne^{\mu^2}\left(\frac{d}{d\mu}\right)^n e^{-\mu^2}=n!\sum_{k=0}^{\lfloor n/2\rfloor}\frac{(-1)^k}{(n-2k)!k!}(2\mu)^{n-2k},\quad n\in\N_0,
\end{equation}
where $\lfloor\cdot\rfloor$ denotes the floor of a number  (i.e. the integral part of a number).

Let $\psi_n(\mu)=\mathcal H_n(\mu) e^{-\mu^2/2}$, $\mu\in\mathbb R$, $n\in\N_0$. It is well known \cite[18.2.1, 18.2.5 (and Table 18.3.1)]{NIST}   that
\begin{equation}
\label{ort}
\int_{-\infty}^\infty\psi_n(\mu)\psi_m(\mu)\,d\mu=\sqrt{\pi}2^nn!\delta_{mn},\quad n\in\N_0,\ m\in\N_0,
\end{equation}
where $\delta_{mn}$ is the Kronecker delta, and
$\{\psi_n\}_{n=0}^\infty$ is an orthogonal basis in $L^2(\mathbb R)$. It
is easy to see that
\begin{equation}
\label{ft}
\widehat{\psi}_n=\mathscr F\psi_n=(-i)^n\psi_n,\qquad n\in\N_0.
\end{equation}
Let $\alpha>0$. Put 
\begin{align}
\label{psiext000}	
\psi_n^\alpha(\mu)&=\frac{1}{\big(\sqrt{2\alpha\pi}2^nn!\big)^{1/2}}
\psi_n\left(\frac{\mu}{\sqrt{2\alpha}}\right),\quad \mu\in\mathbb R,\ n\in\N_0,
\\ 
\label{psiext0} 
\widehat{\psi}_n^\alpha(\lambda)&=\left(\mathscr F\psi_n^\alpha\right)(\lambda),\quad  \lambda\in\mathbb R, \ n\in\N_0.
\end{align}
Then
\begin{equation}
\label{psiext}
\widehat{\psi}_n^\alpha(\lambda)=\left(\frac{2\alpha}{\pi}\right)^{1/4}\frac{(-i)^n}{\big(2^nn!\big)^{1/2}}\psi_n(\sqrt{2\alpha}\lambda),\quad\lambda\in\mathbb R,\ n\in\N_0.
\end{equation}
According to \eqref{ort},  we get
\begin{equation}
\label{ortt}
\Big\langle\psi_n^\alpha,\psi_m^\alpha\Big\rangle=\Big\langle\widehat{\psi}_n^\alpha,\widehat{\psi}_m^\alpha\Big\rangle=\delta_{mn},\quad n\in\N_0,\ m\in\N_0.
\end{equation}
Thus, $\Big\{\psi_n^\alpha\Big\}_{n=0}^\infty$ and
$\Big\{\widehat{\psi}_n^\alpha\Big\}_{n=0}^\infty$ are orthonormal bases in $H^0(\mathbb{R})=H_0(\mathbb{R})=L^2(\mathbb{R})$. In addition, $\Big\{\psi_{2p+1}^\alpha\Big\}_{p=0}^\infty$ and
$\Big\{\widehat{\psi}_{2p+1}^\alpha\Big\}_{p=0}^\infty$ are orthonormal bases in $\widetilde H^0(\R)=\widetilde H_0(\R)$. These bases were considered in \cite{FKh}.

Let $T^*>T$ be fixed. Put 
\begin{equation}
\label{basis2}
\Theta_{nm}^{TT^*}\big((\cdot)_{[1]}, (\cdot)_{[2]}\big)=\psi_{2n+1}^T\big((\cdot)_{[1]}\big) \psi_m^{T*}\big((\cdot)_{[2]}\big), \quad n\in\N_0,\ m\in\N_0.
\end{equation}
Since $\{\psi_p^\alpha\}_{p=0}^\infty$ is an orthonormal basis in $H^0(\R)=L^2(\R)$ (for any $\alpha>0$), we obtain that $\{\Theta_{nm}^{TT^*}\}_{n,m=0}^\infty$ is an orthonormal basis in $\widetilde H^0\left(\R^2\right)\subset L^2\left(\R^2\right)$.

Consider an auxiliary control system for the heat equation on $\R$:
\begin{align}
\label{eqq}
y_t&= y_{x_1 x_1} -2\upsilon\delta', \quad t\in[0,T],
\\
\label{icc}
y(\cdot,0)&=0,
\end{align}  
where $\delta$ is the Dirac distribution with respect to $x_1$, 
$\left(\frac d{dt}\right)^s y: [0,T] \to  \widetilde H^{-2s}(\R)$, $s=0,1$, $\upsilon:[0,T]\to \R$ is a control, $\upsilon\in L^\infty(0,T)$. 

Put $z(\cdot,t)=\mathscr F_{x_1\to\sigma_1} y(\cdot,t)$, $t\in[0,T]$. Then  problem \eqref{eqq}, \eqref{icc} is equivalent to the problem
\begin{align}
\label{eqqf}
z_t&= -\sigma_1^2 z -\sqrt{\frac2\pi}i\sigma_1 \upsilon, \quad \sigma_1\in\R,\ t\in[0,T],
\\
\label{iccf}
z(\cdot,0)&=0,
\end{align} 
where $\left(\frac d{dt}\right)^s z: [0,T] \to  \widetilde H_{-2s}(\R)$, $s=0,1$.
Put
\begin{equation}
\label{sol0}
z(\sigma_1,t)=-\sqrt{\frac2\pi}i\sigma_1\int_0^t e^{-(t-\xi)\sigma_1^2} \upsilon(\xi)\,d\xi,\quad \sigma_1\in\R,\ t\in[0,T].
\end{equation}
Then $z(\cdot,t)\in \widetilde H_0 (\R )$, and it 
is  a unique solution to \eqref{eqqf}, \eqref{iccf}. Let $q$ be a polynomial. 
Put
\begin{align}
\label{contr3-1}
f(\sigma_2, \xi)&=e^{-(T^*-T+\xi)\sigma_2^2} q(\sigma_2),\quad \sigma_2\in\R,\ \xi\in[0,T],
\\
\label{contr3}
\widehat u(\sigma_2,\xi)&=\upsilon(\xi)f(\sigma_2,\xi)
=\upsilon(\xi)e^{-(T^*-T+\xi)\sigma_2^2} q(\sigma_2),\quad \sigma_2\in\R,\ \xi\in[0,T],
\\
\label{sol3}
V(\sigma,t)&=-\sqrt{\frac2\pi}i\sigma_1\int_0^t e^{-(t-\xi)|\sigma|^2} \widehat u(\sigma_2,\xi)\,d\xi,\quad \sigma\in\R^2,\ t\in[0,T].
\end{align}
With regard to \eqref{sol0}--\eqref{sol3}, we have
\begin{align}
\label{sol3-1}
V(\sigma,t)&=\left(-\sqrt{\frac2\pi}i\sigma_1\int_0^t e^{-(t-\xi)\sigma_1^2} \upsilon(\xi)\,d\xi\right)e^{-(T^*-T+t)\sigma_2^2} q(\sigma_2)
\notag\\
&=z(\sigma_1,t)  f(\sigma_2,t),\quad \sigma\in\R^2,\ t\in[0,T].
\end{align}
Since $z(\cdot,t)\in\widetilde H_0\Big(\R\Big)$  and $f(\cdot,t)\in L^2(\R)$, we have $V(\cdot,t)\in\widetilde H_0\left(\R^2\right)$, $t\in[0,T]$. Moreover, 
\begin{equation}
\label{sol3-2}
V(\sigma,T)=z(\sigma_1,T)f(\sigma_2,T)= z(\sigma_1,T) e^{-T^*\sigma_2^2} q(\sigma_2),\quad \sigma\in\R^2.
\end{equation} 
It follows from \eqref{wu} and \eqref{sol3} that
\begin{equation}
\label{sol3-3} 
\mathscr F_{\sigma\to x}^{-1}V(\cdot,t) =\mathcal{W}_{u}(\cdot,t),\quad t\in[0,T], 
\end{equation} 
where $u(\cdot,t)=\mathscr F_{\sigma_2\to x_2}^{-1}\widehat u(\cdot,t)$ and $u\in U[0,T]$. 

\begin{theorem}
	\label{thclos}
	Let $T^*>T>0$, $n\in\N_0$, and $m\in\N_0$. Then $\Theta_{nm}^{TT^*}\in\overline{\mathscr{R}_T(0)}$.
\end{theorem}

\begin{proof}
	It follows from \cite[Theorem 7.2]{FKh} that there exists a sequence $\big\{\upsilon^n_l\big\}_{l=1}^\infty\subset L^\infty(0,T)$ such that the solutions $y^n_l$ to \eqref{eqq}, \eqref{icc} with  $\upsilon=\upsilon^n_l$ satisfy the condition
	\begin{equation}
	\label{est0}
	\left\|y^n_l(\cdot,T)-\psi_{2n+1}^T\right\|^0 \to 0 \quad \text{as}\ l\to\infty.
	\end{equation} 
	Setting $z^n_l(\cdot,t)=\mathscr F_{x_1\to\sigma_1} y^n_l(\cdot,t)$, $t\in[0,T]$, we conclude that $z^n_l(\cdot,t)\in \widetilde H_0(\R)$, $t\in[0,T]$, and it is the unique solution to \eqref{eqqf}, \eqref{iccf} with  $\upsilon=\upsilon^n_l$, $l\in\N$. Taking into account \eqref{psiext0} and \eqref{est0}, we get 
	\begin{equation}
	\label{est1}
	\left\|z^n_l(\cdot,T)-\widehat\psi_{2n+1}^T\right\|_0 \to 0 \quad \text{as}\ l\to\infty.
	\end{equation}
	Put
	\begin{align}
	\label{contr30}
	f_m(\sigma_2, \xi)&=e^{(T-\xi)\sigma_2^2} \widehat\psi_m^{T^*}(\sigma_2),\quad \sigma_2\in\R,\ \xi\in[0,T],
	\\
	\label{contr31}
	\widehat u_{nm}^l(\sigma_2,\xi)&=\upsilon^n_l(\xi)f_m(\sigma_2, \xi),\quad \sigma_2\in\R,\ \xi\in[0,T],\ l\in\N,
	\\
	\label{sol31}
	V_{nm}^l(\sigma,t)&=-\sqrt{\frac2\pi}i\sigma_1\int_0^t e^{-(t-\xi)|\sigma|^2} \widehat u_{nm}^l(\sigma_2,\xi)\,d\xi,\  \sigma\in\R^2,\ t\in[0,T],\ l\in\N.
	\end{align}
	Note that $f_m(\cdot, \xi)$ is a polynomial multiplied by $e^{-(T^*-T+\xi)(\cdot)^2}$, $\xi\in[0,T]$. With regard to  \eqref{contr3}--\eqref{sol3-1}, we have $V_{nm}^l(\cdot,t)\in\widetilde H_0(\R^2)$ and
	\begin{align}
	\label{sol3_1}
	V_{nm}^l(\sigma,t)=z^n_l(\sigma_1,t) f_m(\sigma_2, t),\quad \sigma\in\R^2,\ t\in[0,T],\ l\in\N.
	\end{align}
	According to \eqref{ortt} and \eqref{est1},  we get
	\begin{align}
	\label{est2}
	\left\|\widehat \Theta_{nm}^{TT*}-V_{nm}^l(\cdot,T)\right\|_0 
	&= \left\|z^n_l(\cdot,T)-\widehat\psi_{2n+1}^T\right\|_0\left\|\widehat\psi_m^{T^*}\right\|_0
	\notag\\
	&=\left\|z^n_l(\cdot,T)-\widehat\psi_{2n+1}^T\right\|_0 \to 0 \quad \text{as}\ l\to\infty,
	\end{align}
	where
	\begin{align}
	\label{basis3}
	\widehat \Theta_{nm}^{TT^*}= \mathscr F \Theta_{nm}^{TT^*}
	&=\left(\mathscr F_{x_1\to\sigma_1} \psi_{2n+1}^T\right)((\cdot)_{[1]}) \left(\mathscr F_{x_2\to\sigma_2} \psi_m^{T^*}\right)((\cdot)_{[2]})
	\notag\\
	&=\widehat \psi_{2n+1}^T((\cdot)_{[1]}) \widehat \psi_m^{T^*}((\cdot)_{[2]}).
	\end{align}
	Here  \eqref{psiext0} and \eqref{basis2} are taken into account. 
	Due to \eqref{sol3-3}, we have $\mathscr F_{\sigma\to x}^{-1}V_{nm}^l(\cdot,t) =\mathcal{W}_{u_{nm}^l}(\cdot,t)$, $t\in[0,T]$, where $u_{nm}^l(\cdot,t)=\mathscr F_{\sigma_2\to x_2}^{-1}\widehat u_{nm}^l(\cdot,t)$, $l\in\N$. Therefore, \eqref{est2} yields 
	\begin{equation}
	\label{est4}
	\left\| \Theta_{nm}^{TT*}-\mathcal{W}_{u_{nm}^l}(\cdot,T) \right\|^0\to 0  \quad \text{as}\ l\to\infty.
	\end{equation}
	
	Finally, let us calculate $u_{nm}^l$ to ensure $u_{nm}^l\in U[0,T]$, $l\in\N$. Let $l\in\N$ be fixed. Put $\Delta T= T^*-T$.
	We have 
	\begin{align}
	\label{lll2}
	&\mathscr F_{\sigma_2\to x_2}^{-1} \left(e^{(T-\xi)(\cdot)^2} \widehat\psi_m^{T^*}\right)(x_2)
	\notag\\ 
	&\quad
	=\left(\frac{2T^*}{\pi}\right)^{1/4} \frac{(-i)^m}{(2^m m!)^{1/2}} \mathscr F_{\sigma_2\to x_2}^{-1} \left(\mathcal H_m\Big(\sqrt{2T^*}(\cdot)\Big) e^{-(\Delta T+\xi)(\cdot)^2}\right)(x_2)
	\notag\\ 
	&\quad
	= \left(\frac{2T^*}{\pi}\right)^{1/4} \frac{(-i)^m}{(2^m m!)^{1/2}}\mathcal H_m\left(-i\sqrt{2T^*}\frac{\partial}{\partial x_2}\right) \frac{e^{-x_2^2/(4(\Delta T+\xi))}}{\sqrt{2(\Delta T+\xi)}}
	\notag\\ 
	&\quad
	= \left(\frac{2T^*}{\pi}\right)^{1/4} \frac{(-1)^m m!}{(2^m m!)^{1/2}} \frac{1}{\sqrt{2(\Delta T+\xi)}}
	\notag\\ 
	&\qquad
	\times 
	\sum_{k=0}^{\lfloor m/2\rfloor}
	\frac{(2\sqrt{2T^*})^{m-2k}}{(m-2k)! k!}
	\left(\frac{\partial}{\partial x_2}\right)^{m-2k} e^{-x_2^2/(4(\Delta T+\xi))},
	\notag\\
	&\kern50ex 
	x_2\in\R,\ \xi\in[0,T].
	\end{align}
	From \eqref{hermite},  it follows that
	\begin{equation}
	\label{lll3}
	\left(\frac{d}{d\mu}\right)^p e^{-\mu^2}= (-1)^p \mathcal H_p(\mu)e^{-\mu^2},\quad \mu\in\R,\ p\in\N_0.
	\end{equation}
	Therefore,
	\begin{align}
	\label{lll4}
	\left(\frac{\partial}{\partial x_2}\right)^p e^{-x_2^2/(4(\Delta T+\xi))}= \frac{(-1)^p}{(2\sqrt{\Delta T+\xi})^p} \mathcal H_p&\left(\frac{x_2}{2\sqrt{\Delta T+\xi}}\right)e^{-x_2^2/(4(\Delta T+\xi))},
	\notag\\
	& x_2\in\R,\ \xi\in[0,T],\ p\in\N_0.
	\end{align}
	With regard to \eqref{lll2}, we obtain
	\begin{align*}
	&\mathscr F_{\sigma_2\to x_2}^{-1} \left(e^{(T-\xi)(\cdot)^2} \widehat\psi_m^{T^*}\right)(x_2)
	\\
	&\quad
	= \left(\frac{2T^*}{\pi}\right)^{1/4} \left(\frac{ m!}{2^m }\right)^{1/2} \frac{1}{\sqrt{2(\Delta T+\xi)}}
	\\
	&\qquad\times
	\sum_{k=0}^{\lfloor m/2\rfloor}
	\frac{e^{-x_2^2/(4(\Delta T+\xi))}}{(m-2k)! k!}
	\mathcal H_{m-2k}\left(\frac{x_2}{2\sqrt{\Delta T+\xi}}\right)
	\left(\frac{2T^*}{\Delta T+\xi}\right)^{m/2-k} 
	\\
	&\quad
	= \left(\frac{2T^*}{\pi}\right)^{1/4} \frac1{(2^m m!)^{1/2}} \frac{e^{-x_2^2/(4(\Delta T+\xi))}}{\sqrt{2(\Delta T+\xi)}}\left(\frac{T^*+T-\xi}{\Delta T+\xi}\right)^{m/2}
	\left(\frac{2T^*}{T^*+T-\xi}\right)^{m/2}
	\\
	&\qquad
	\times
	\sum_{k=0}^{\lfloor m/2\rfloor}
	\left(1-\frac{T^*+T-\xi}{2T^*}\right)^k
	\frac{m!}{k!(m-2k)!}
	\mathcal H_{m-2k}\left(\frac{x_2}{2\sqrt{\Delta T+\xi}}\right),
	\\
	&\kern60ex
	 x_2\in\R,\ \xi\in[0,T].
	\end{align*}
	Taking into account the multiplication theorem for Hermite polynomials (see, e.g. \cite[18.18.13]{NIST}):
	$$
	\mathcal H_m(\lambda\mu)
	=\lambda^m \sum_{k=0}^{\lfloor m/2\rfloor} \frac{m!}{k!(m-2k)!} \left(1-\frac{1}{\lambda^2}\right)^k \mathcal H_{m-2k}(\mu),\quad \mu\in\R,\ \lambda\neq0,
	$$
	and setting $\lambda=\sqrt{\frac{2T^*}{T^*+T-\xi}}$, $\mu=\frac{x_2}{2\sqrt{\Delta T+\xi}}$, we get
	\begin{align}
	\label{lll5a}
	\mathscr F_{\sigma_2\to x_2}^{-1} & \left(e^{(T-\xi)(\cdot)^2} \widehat\psi_m^{T^*}\right)(x_2)
	\notag\\
	&=\left(\frac{2T^*}{\pi}\right)^{1/4} \frac1{(2^m m!)^{1/2}} \frac{e^{-x_2^2/(4(\Delta T+\xi))}}{\sqrt{2(\Delta T+\xi)}}\left(\frac{T^*+T-\xi}{\Delta T+\xi}\right)^{m/2}
	\notag \\ 
	&\quad
	\times
	\mathcal
	H_m\left(x_2\sqrt{\frac{T^*}{2(\Delta T+\xi)(T^*+T-\xi)}}\right),
	\quad  x_2\in\R,\ \xi\in[0,T].
	\end{align}
	Due to \eqref{contr31}, we have
	\begin{align}
	\label{lll5}
	u_{nm}^l(x_2,\xi)&=\mathscr F_{\sigma_2\to x_2}^{-1} \left(e^{(T-\xi)(\cdot)^2} \widehat\psi_m^{T^*}\right)(x_2) \upsilon_l^n(\xi)
	\notag\\ 
	&
	=\left(\frac{2T^*}{\pi}\right)^{1/4} \frac{\upsilon_l^n(\xi)}{(2^m m!)^{1/2}} \frac{e^{-x_2^2/(4(\Delta T+\xi))}}{\sqrt{2(\Delta T+\xi)}}\left(\frac{T^*+T-\xi}{\Delta T+\xi}\right)^{m/2}
	\notag\\ 
	&
	\quad\times
	\mathcal H_m\left(x_2\sqrt{\frac{T^*}{2(\Delta T+\xi)(T^*+T-\xi)}}\right),
	\ 
	x_2\in\R,\ \xi\in[0,T].
	\end{align}
	
	Evidently, $u_{nm}^l\in U[0,T]$, $l\in\N$.
	The theorem is proved.
\end{proof} 

\begin{proof}[Proof of Theorem \ref{thappc}]
	With regard  to Definition \ref{def-appr} and Theorems \ref{th-sol}, \ref{reachprop}, it is sufficient to show that $\widetilde{H}^0\left(\mathbb{R}^2\right)=\overline{\mathscr{R}_T(0)}$ to prove the theorem. Set $T^*>T$. Since $\{\Theta_{nm}^{TT^*}\}_{n,m=0}^\infty$ is an orthonormal basis in $\widetilde H^0\left(\R^2\right)$, each $f\in \widetilde{H}^0\left(\mathbb{R}^2\right)$ can be approximated in this space by the sums
	$$
	\sum_{n=0}^N \sum_{m=0}^M f_{nm}\Theta_{nm}^{TT^*},
	$$
	where $N\in\N$, $M\in\N$. With regard to Theorem \ref{thclos}, we conclude that $f\in \overline{\mathscr{R}_T(0)}$ that was to be proved.
\end{proof}


\section{Numerical solution to the approximate controllability problem}
\label{sect5}

In this section, we construct controls approximately steering to a state $W^T\in\widetilde{H}^0\left(\mathbb{R}^2\right)$ from a state $W^0\in\widetilde{H}^0\left(\mathbb{R}^2\right)$ by using the proof of Theorem \ref{thclos}. 

Put $W_0^T=W^T-\mathcal W_0(\cdot,T)$, $V_0^T=\mathscr F W_0^T$. According to Theorem \ref{th-sol}, $W_0^T\in\widetilde H^0\left(\R^2\right)$, hence $V_0^T\in\widetilde H_0\left(\R^2\right)$. Due to Theorem \ref{reachprop}, we have to construct controls $\{u_k\}_{k=1}^\infty\subset U[0,T]$ such that
\begin{equation}
\label{num0}
\left\|W_0^T- \mathcal W_{u_k}(\cdot,T)\right\|^0\to 0 \quad \text{as}\ k\to\infty.
\end{equation}
Set $T^*>T$.
Since $\big\{\Theta_{nm}^{TT^*}\big\}_{n,m=0}^\infty$ is an orthonormal basis in $\widetilde H^0\left(\R^2\right)$, we have
\begin{equation}
\label{num3}
\left\|W_0^T-W_{NM} \right\|^0\to0 \quad\text{as}\ (N,M)\to\infty
\end{equation}
for
\begin{equation}
\label{num2}
W_{NM}=\sum_{n=0}^N \sum_{m=0}^M W_{nm}\Theta_{nm}^{TT^*},\quad N\in\N,\ M\in \N,
\end{equation}
where $W_{nm}=\left\langle W_0^T,\Theta_{nm}^{TT^*} \right\rangle$, $n\in\N_0$, $m\in\N_0$.

Let an arbitrary $\varepsilon>0$ be fixed. Then there exists $N\in\N$ and $M\in\N$ such that 
\begin{equation}
\label{num4}
\left\|W_0^T-W_{NM} \right\|^0<\varepsilon/2.
\end{equation}

Let us construct controls $\{u_k\}_{k=1}^\infty\subset U[0,T]$ for which 
\begin{equation}
\label{num5}
\left\|W_{NM}-\mathcal W_{u_k}(\cdot,T) \right\|^0\to0 \quad \text{as}\ k\to \infty.
\end{equation}
To this aid, we use the method proposed in \cite{FKh} to construct controls for the heat equation oh a half-axis.
With regard to \eqref{basis3} and \eqref{num2}, we have 
\begin{align}
\label{num6}
W_{NM}&=\sum_{n=0}^N \sum_{m=0}^M W_{nm}\psi_{2n+1}^T\big((\cdot)_{[1]}\big) \psi_m^{T^*}\big((\cdot)_{[2]}\big)
\notag\\
&=\sum_{m=0}^M\omega_m^N\big((\cdot)_{[1]}\big)\psi_m^{T^*}\big((\cdot)_{[2]}\big),\quad N\in\N,\ M\in \N,
\end{align}
where
\begin{equation}
\label{num7}
\omega_m^N
=\sum_{n=0}^N W_{nm}\psi_{2n+1}^T, \quad m=\overline{0,M}.
\end{equation}
Put
\begin{align}
&\varphi_{2p+1}(\lambda)=i\lambda^{2p+1}e^{-T\lambda^2}, \quad \lambda\in\mathbb{R},\  p\in\N_0,\label{fi}
\\
&\varphi_{2p+1}^l(\lambda)=i\lambda^{2p+1} e^{-T\lambda^2}\left(\frac{e^{\lambda^2/l}-1}{\lambda^2/l}\right)^{p+1}, \quad
\lambda\in\R,\ p\in\N_0, \ l\in\N, 
\nonumber\\
&v_l^p(\xi)=\begin{cases}
\displaystyle
(-1)^{p-j}\binom{p}{j} l^{p+1},& \displaystyle \xi\in\left(\frac{j}{l},\frac{j+1}{l}\right),\
j=\overline{0,p}\\
0,& \displaystyle \xi\notin\left( 0,\frac{p+1}{l}\right)
\end{cases},\quad  p\in\N_0, \ l\in\N.
\label{contrl}
\end{align}
Note that $v_l^p\to(-1)^p\delta^{(p)}$ as $l\to\infty$ in $H^{-1}(\R)$ for  each $p=\overline{0,\infty}$.
Due to \cite{FKh}, we get
\begin{equation}
\label{aphi}
\forall p\in\N_0\quad
\left\|\varphi_{2p+1}^{}-\varphi_{2p+1}^l\right\|_0\to 0 \quad\text{as}\ l\to\infty.
\end{equation}
Let $p\in\N_0$, $l\in\N$, and let $\mathbf y_l^p$  be the solution to \eqref{eqq}, \eqref{icc} with $\upsilon=v_l^p$. Put $\mathbf z_l^p(\cdot,t)=\mathscr F \mathbf y_l^p(\cdot,t)$, $t\in[0,T]$. Then $\mathbf z_l^p$ is the solution to \eqref{eqqf}, \eqref{iccf}, and  \eqref{sol0} holds.
It is easy to see that
\begin{equation}
\label{zlnfi}
\mathbf z_l^p(\sigma_1,T) =-\sqrt{\frac2\pi}\varphi_{2p+1}^l(\sigma_1),\quad \sigma_1\in\R.
\end{equation}
Taking into account \eqref{aphi}, we get
\begin{equation*}
\left\|\mathbf z_l^p(\cdot,T)+\sqrt{\frac2\pi}\varphi_{2p+1}\right\|_0\to 0 \quad\text{as}\ l\to\infty.
\end{equation*}

Put $V_{NM} =\mathscr F_{x\to\sigma} W_{NM} $. It follows from \eqref{num6} and \eqref{num7} that 
\begin{equation}
\label{vtnm}
V_{NM} =\sum_{m=0}^M\widehat\omega_m^N\big((\cdot)_{[1]}\big)\widehat\psi_m^{T^*}\big((\cdot)_{[2]}\big),\quad N\in\N,\ M\in \N,
\end{equation}
where 
\begin{equation}
\label{omegan}
\widehat\omega_m^N
=\sum_{n=0}^N W_{nm}\widehat\psi_{2n+1}^T, \quad m=\overline{0,M}.
\end{equation}
Using \eqref{hermite}, \eqref{psiext0}, \eqref{psiext}, and \eqref{fi}, we obtain 
\begin{align*}
\widehat\omega_m^N(\sigma_1)
&=\sum_{n=0}^N W_{nm}(-i)^{2n+1}\left(\frac{2T}{\pi}\right)^{1/4}\frac{e^{-T\sigma_1^2}\mathcal H_{2n+1}(\sqrt{2T}\sigma_1)}{\sqrt{2^{2n+1}(2n+1)!}}\\
&=-i\left(\frac{2T}{\pi}\right)^{1/4}e^{-T\sigma_1^2}
\\
&\quad\times
\sum_{n=0}^N W_{nm}(-1)^n\sqrt{\frac{(2n+1)!}{2^{2n+1}}}\sum_{j=0}^n\frac{(-1)^j(2\sqrt{2T}\sigma_1)^{2n+1-2j}}{j!(2n+1-2j)!}\\
&=-\left(\frac{2T}{\pi}\right)^{1/4}\sum_{n=0}^N W_{nm}\sqrt{\frac{(2n+1)!}{2^{2n+1}}}\sum_{p=0}^n\frac{(-1)^p(2\sqrt{2T})^{2p+1}}{(n-p)!(2p+1)!}i\sigma_1^{2p+1}e^{-T\sigma_1^2}\\
&=\sum_{p=0}^N\varphi_{2p+1}(\sigma_1)\sum_{n=p}^N W_{nm}h_p^n
=\sum_{p=0}^N g_{pm}^N\varphi_{2p+1}(\sigma_1), \quad \sigma_1\in\mathbb{R},\  m=\overline{0,M},
\end{align*} 
where
\begin{align}
\label{gpm}
g_{pm}^N&=\sum_{n=p}^N W_{nm}h_p^n,&& p=\overline{0,N},\ m=\overline{0,M},
\\
\label{hpn}
h_p^n&=\left(\frac{2T}{\pi}\right)^{1/4}\frac{(-1)^{p+1}(2\sqrt{2T})^{2p+1}}{(n-p)!(2p+1)!}\sqrt{\frac{(2n+1)!}{2^{2n+1}}},&& p=\overline{0,N},\ n=\overline{p,N}.
\end{align}
Thus,
\begin{equation}
\label{vtnmser}
V_{NM} =\sum_{m=0}^M\sum_{p=0}^N g_{pm}^N\varphi_{2p+1}\big((\cdot)_{[1]}\big) \widehat\psi_m^{T^*}\big((\cdot)_{[2]}\big).
\end{equation}
Put
\begin{equation}
\label{vtnmlser}
V_{NM}^{l}=\sum_{m=0}^M\sum_{p=0}^N g_{pm}^N\varphi_{2p+1}^l\big((\cdot)_{[1]}\big) \widehat\psi_m^{T^*}\big((\cdot)_{[2]}\big),\quad l\in\mathbb{N}.
\end{equation}
Taking into account \eqref{aphi}, we see that there exists  $l\in \N$ such that 
\begin{equation}
\label{normv}
\left\|V_{NM}-V_{NM}^{l} \right\|_0\leq\sum_{m=0}^M\sum_{p=0}^N \left|g_{pm}^N\right|\left\|\varphi_{2p+1}-\varphi_{2p+1}^l\right\|_0 \left\|\widehat\psi_m^{T^*}\right\|_0<\varepsilon/2.
\end{equation}
Using  \eqref{zlnfi}, we get
\begin{align}
\label{vnml}
V_{NM}^{l}(\sigma)&=-\sqrt{\frac{\pi}{2}}
\sum_{m=0}^M\sum_{p=0}^N g_{pm}^N \mathbf z_l^p(\sigma_1,T)\widehat\psi_m^{T^*}(\sigma_2)
,\quad \sigma\in\mathbb{R}^2. 
\end{align}
Put 
\begin{align}
\label{hatumnl}
\widehat{u}_{NM}^l(\cdot,\xi)&=-\sqrt{\frac{\pi}{2}}\sum_{m=0}^M\sum_{p=0}^N g_{pm}^N v_l^p(\xi) f_m(\cdot,\xi)
,\quad  \xi\in[0,T],
\\
\label{contr4}
u_{NM}^l(\cdot,\xi)&=\mathscr F^{-1}_{\sigma_2\to x_2}\widehat{u}_{NM}^l(\cdot,\xi),\quad  \xi\in[0,T],
\end{align}
where $f_m$ is defined by \eqref{contr30}. 
With regard to \eqref{contr3-1}, \eqref{contr3}, \eqref{sol3-2}, \eqref{sol3-3}, and \eqref{contr30}, we conclude that 
\begin{equation}
\label{wnml}
\mathscr F^{-1} V_{NM}^{l}=\mathcal W_{u_{NM}^l}(\cdot,T).
\end{equation}
Taking into account \eqref{normv}, we have
\begin{equation}
\label{normw}
\left\|W_{NM}-\mathcal W_{u_{NM}^l}(\cdot,T)\right\|^0
=\left\|V_{NM}-V_{NM}^{l}(\cdot,T)\right\|_0<\varepsilon/2.
\end{equation}
With regard to \eqref{lll5a} and \eqref{contr4}, we obtain
\begin{align}
\label{umnl}
u_{NM}^l&(x_2,\xi)
\notag\\
&
=-\left(\frac{\pi T^*}{2}\right)^{1/4}\sum_{m=0}^M\sum_{p=0}^N g_{pm}^N  \frac{v_l^p(\xi)}{(2^m m!)^{1/2}} \frac{e^{-x_2^2/(4(\Delta T+\xi))}}{\sqrt{2(\Delta T+\xi)}}\left(\frac{T^*+T-\xi}{\Delta T+\xi}\right)^{m/2}
\notag\\ 
&\quad
\times
\mathcal H_m\left(x_2\sqrt{\frac{T^*}{2(\Delta T+\xi)(T^*+T-\xi)}}\right),
\quad x_2\in\mathbb{R},\ \xi\in[0,T],
\end{align} 
where $\Delta T=T^*-T$, $g_{pm}^N$ is defined by \eqref{gpm}, $v_l^p$ is defined by \eqref{contrl}, $p=\overline{0,N}$, $m=\overline{0,M}$, $l\in\mathbb{N}$, $N\in\N$, $M\in \N$.
Due to \eqref{num4} and \eqref{normw}, we conclude that
\begin{equation}
\label{lll9}
\left\|W^T-\left(\mathcal W_0(\cdot,T)+\mathcal W_{u_{NM}^{l}}(\cdot,T)\right) \right\|^0
=\left\|W_0^T-\mathcal W_{u_{NM}^{l}}(\cdot,T)\right\|^0<\varepsilon,
\end{equation}
therefore the controls $u_{NM}^{l}$, $N\in\N$, $M\in\N$, $l\in\N$, solves the approximate controllability problem for system \eqref{eq1}, \eqref{ic1}.

Note that we have found the approximate end state  in the form
\begin{equation}
\label{apprw}
\mathcal W_0(\cdot,T)+\mathcal W_{u_{NM}^{l}}(\cdot,T), \quad l\in\mathbb{N},\ N\in\N,\ M\in \N,
\end{equation}
and the control $u_{NM}^l$ solving the approximate controllability problem for considered system in the form \eqref{umnl}.


\section{Example}
\label{ex}

\begin{example}
	\label{ex3}
	Let $T=2$, $x_1>0$, $x_2\in\mathbb{R}$, 
	\begin{align*}
	&w^0(x)=-\frac{9}{8}e^{-5/8}\left(\frac{1}{8\pi T^3}\right)^{1/4}x_1e^{-|x|^2/(8T)},\\
	&w^T(x)=-e^{-5/8}\left(\!\frac{1}{8\pi T^3}\right)^{1/4}\! \left(e^{-|x|^2/(4T)}e^{-x_2/\sqrt{T}}\sinh\left(\frac{x_1}{2\sqrt{T}}\right)+\frac{x_1}{2}e^{-|x|^2/(12T)}\!\right).
	\end{align*} 
	Let $W^0$ and $W^T$ be the odd extensions of $w^0$ and $w^T$  with respect to $x_1$.
	Consider system \eqref{eq1}, \eqref{ic1}. We can see that $W^0\in\widetilde{H}^0\left(\mathbb{R}^2\right)$ and $W^T\in\widetilde{H}^0\left(\mathbb{R}^2\right)$. Due to Theorem \ref{thappc}, we conclude that the state $W^0$ is approximately controllable to the state $W^T$ in the time $T=2$. Let us find  controls $u_{NM}^l$  approximately steering the state $W^0$ to the state $W^T$ and the approximate end states associated with these controls.
	
	Put $W^T_0=W^T-\mathcal W_0(\cdot,T)$ (see \eqref{wo}), $V^T_0=\mathscr{F}
	_{x\rightarrow\sigma} W^T_0$. One can easily obtain 
	\begin{align} 
	&\mathcal W_0(x,T)=-e^{-5/8}\left(\frac{1}{8\pi T^3}\right)^{1/4}\frac{x_1}{2}e^{-|x|^2/(12T)},&& x\in\mathbb{R}^2,\label{w0ex}\\
	&W^T_0(x)=-e^{-5/8}\left(\frac{1}{8\pi T^3}\right)^{1/4}e^{-|x|^2/(4T)}e^{-x_2/\sqrt{T}}\sinh\left(\frac{x_1}{2\sqrt{T}}\right), && x\in\mathbb{R}^2.
	\notag
	\end{align}
	Consider any $T^*>T$. Let $N\in\N$ and $M\in\N$ be fixed, and let $W_{NM}$ be determined by \eqref{num2}, 
	where $W_{nm}=\left\langle W_0^T,\Theta_{nm}^{TT^*} \right\rangle$, $n\in\N_0$, $m\in\N_0$. Let us find the coefficients. We have
	\begin{align}
	\label{wnmcoef00} 
	W_{nm}&=-e^{-5/8}\left(\frac{1}{8\pi T^3}\right)^{1/4}\int_{-\infty}^\infty e^{-x_1^2/(4T)}\sinh\left(\frac{x_1}{2\sqrt{T}}\right)\psi_{2n+1}^T(x_1)\,dx_1
	\notag\\
	&\quad\times 
	\int_{-\infty}^\infty e^{-x_2^2/(4T)}e^{-x_2/\sqrt{T}}\psi_m^{T^*}(x_2)\,dx_2,\quad n\in\N_0,\ m\in\N_0.
	\end{align}
	Using \eqref{hermite} and \eqref{psiext000}, we obtain
	\begin{align} 
	&\int_{-\infty}^\infty e^{-x_1^2/(4T)}\sinh\left(\frac{x_1}{2\sqrt{T}}\right)\psi_{2n+1}^T(x_1)\,dx_1
	\nonumber\\
	&=\left(\frac{1}{2\pi T}\right)^{1/4}\sqrt{\frac{(2n+1)!}{2^{2n+1}}}\int_{-\infty}^\infty e^{-x_1^2/(2T)}\sinh\left(\frac{x_1}{2\sqrt{T}}\right)
	\notag\\
	&\quad\times
	\sum_{k=0}^n
	\frac{(-1)^k}{(2n-2k+1)!k!}\left(\sqrt{\frac{2}{T}}x_1\right)^{2n-2k+1}\,dx_1\nonumber\\
	&=\left(\frac{1}{2\pi T}\right)^{1/4}\sqrt{\frac{(2n+1)!}{2^{2n+1}}}\sqrt{\frac{2}{T}}\sum_{p=0}^n \frac{(-1)^{n-p}}{(2p+1)!(n-p)!}\left(\frac{2}{T}\right)^p
	\notag\\
	&\quad\times
	\int_{-\infty}^\infty x_1^{2p+1}e^{-x_1^2/(2T)}\sinh\left(\frac{x_1}{2\sqrt{T}}\right)\,dx_1,\quad n\in\N_0.\label{int11}
	\end{align}  
	Taking into account the integral representation (\cite[18.10.10]{NIST}): 
	\begin{equation}
	\label{herm1}
	\mathcal H_n(x)=\frac{(-2i)^n e^{x^2}}{\sqrt{\pi}}\int_{-\infty}^\infty e^{-t^2}t^n e^{2ixt}\,dt,\quad x\in\mathbb{R},\ n\in\N_0,
	\end{equation}
	for $n=2p+1$, we get
	\begin{align*} 
	\int_{-\infty}^\infty x_1^{2p+1}e^{-x_1^2/(2T)}&\sinh\left(\frac{x_1}{2\sqrt{T}}\right)\,dx_1
	\\
	&=i\sqrt{\pi}e^{1/8}\frac{(-1)^{p+1}T^{p+1}}{2^p}\mathcal H_{2p+1}\left(\frac{i}{2\sqrt{2}}\right),\quad p=\overline{0,n},\ n\in\N_0.
	\end{align*}
	Continuing \eqref{int11}, we have
	\begin{align*} 
	\int_{-\infty}^\infty e^{-x_1^2/(4T)}&\sinh\left(\frac{x_1}{2\sqrt{T}} \right)\psi_{2n+1}^T(x_1)\,dx_1
	\\
	&=\left(\frac{1}{2\pi T}\right)^{1/4}\sqrt{2\pi T}ie^{1/8}(-1)^{n+1}\sqrt{\frac{(2n+1)!}{2^{2n+1}}}
	\\
	&\quad\times
	\sum_{p=0}^n \frac{1}{(2p+1)!(n-p)!}\mathcal H_{2p+1}\left(\frac{i}{2\sqrt{2}}\right)
	\\
	&=(2\pi T)^{1/4}ie^{1/8}(-1)^{n+1}\sqrt{\frac{(2n+1)!}{2^{2n+1}}}
	\\
	&\quad\times
	\sum_{k=0}^n \frac{1}{k!(2n+1-2k)!}\mathcal H_{2n+1-2k}\left(\frac{i}{2\sqrt{2}}\right),\quad n\in\N_0.
	\end{align*}
	Using the connection formula (\cite[18.18.20]{NIST})
	\begin{equation*}
	(2x)^n=\sum_{l=0}^{\lfloor n/2\rfloor}\frac{n!}{l!(n-2l)!}\mathcal H_{n-2l}(x),\quad x\in\mathbb{R},\ n\in\N_0,
	\end{equation*}
	we get 
	\begin{align} 
	\label{int3}
	\int_{-\infty}^\infty e^{-x_1^2/(4T)}&\sinh\left(\frac{x_1}{2\sqrt{T}}\right)\psi_{2n+1}^T(x_1)\,dx_1
	\notag\\
	&=(2\pi T)^{1/4}ie^{1/8}(-1)^{n+1}\sqrt{\frac{(2n+1)!}{2^{2n+1}}}\frac{2^{2n+1}}{(2n+1)!}\left(\frac{i}{2\sqrt{2}}\right)^{2n+1}
	\nonumber\\
	&=\frac{(2\pi T)^{1/4}e^{1/8}}{2^{2n+1}\sqrt{(2n+1)!}},\quad n\in\N_0.
	\end{align}
	
	Using again \eqref{hermite} and \eqref{psiext000}, we obtain
	\begin{align}
	\label{int4} 
	\int_{-\infty}^\infty e^{-x_2^2/(4T)}&e^{-x_2/\sqrt{T}}\psi_m^{T^*}(x_2)\,dx_2
	\nonumber\\
	&=\left(\frac{1}{2\pi T^*}\right)^{1/4}\sqrt{\frac{m!}{2^m}}\sum_{k=0}^{\lfloor m/2\rfloor}\frac{(-1)^k}{k!(m-2k)!}\left(\sqrt{\frac{2}{T^*}}\right)^{m-2k}
	\notag\\
	&\quad\times
	\int_{-\infty}^\infty x_2^{m-2k} e^{-x_2^2(T^*+T)/(4TT^*)-x_2/\sqrt{T}}\,dx_2,\quad m\in\N_0.
	\end{align}
	Taking again into account \eqref{herm1}, we have
	\begin{align*} 
	\int_{-\infty}^\infty x_2^{m-2k} &e^{-x_2^2(T^*+T)/(4TT^*)-x_2/\sqrt{T}}\,dx_2
	\\
	&
	=e^{T^*/(T^*+T)}\left(2\sqrt{\frac{TT^*}{T^*+T}}\right)^{m-2k+1}
	\!\!
	\frac{\sqrt{\pi}}{(2i)^{m-2k}}
	\mathcal H_{m-2k}\left(-i\sqrt{\frac{T^*}{T^*+T}}\right)
	\\
	&
	=2\sqrt{\pi}e^{T^*/(T^*+T)}\frac{(-1)^k}{i^m}\left(\sqrt{\frac{TT^*}{T^*+T}}\right)^{m-2k+1}
	\!\!
	\mathcal H_{m-2k}\left(-i\sqrt{\frac{T^*}{T^*+T}}\right)\!,
	\\
	& \kern47ex 
	k=\overline{0,\lfloor m/2\rfloor},\ m\in\N_0.
	\end{align*} 
	Continuing \eqref{int4}, we get
	\begin{align} 
	&\int_{-\infty}^\infty e^{-x_2^2/(4T)}e^{-x_2/\sqrt{T}}\psi_m^{T^*}(x_2)\,dx_2\nonumber\\
	&=(8\pi T^*)^{1/4}e^{T^*/(T^*+T)}\sqrt{\frac{T}{T^*+T}}i^{-m}\sqrt{\frac{m!}{2^m}}\sum_{k=0}^{\lfloor m/2\rfloor}\frac{1}{k!(m-2k)!}\left(\sqrt{\frac{2T}{T^*+T}}\right)^{m-2k}
	\notag\\
	&\quad\times
	\mathcal H_{m-2k}\left(-i\sqrt{\frac{T^*}{T^*+T}}\right),\quad m\in\N_0.\label{int5}
	\end{align}
	Using the multiplication theorem (\cite[18.18.13]{NIST})  
	\begin{equation*}
	\mathcal H_n(\lambda\mu)=\lambda^n\sum_{j=0}^{\lfloor n/2\rfloor}\frac{n!}{j!(n-2j)!}(1-\lambda^{-2})^j\mathcal H_{n-2j}(\mu),\quad \lambda\in\mathbb{R},\ \mu\in\mathbb{R},\ n\in\N_0,
	\end{equation*}
	and setting $\lambda=i\sqrt{\frac{2T}{\Delta T}}$ and $\mu=-i\sqrt{\frac{T^*}{T^*+T}}$, we obtain
	\begin{align*} 
	\sum_{k=0}^{\lfloor m/2\rfloor}&\frac{1}{k!(m-2k)!} \left(\sqrt{\frac{2T}{T^*+T}}\right)^{m-2k} \mathcal H_{m-2k}\left(-i\sqrt{\frac{T^*}{T^*+T}} \right)
	\\
	&=\left(\frac{2T}{T^*+T}\right)^{m/2}\frac{1}{m!i^m} \left(\frac{\Delta T}{2T}\right)^{m/2}\left(i\sqrt{\frac{2T}{\Delta T}}\right)^m
	\\
	&\quad\times
	\sum_{k=0}^{\lfloor m/2\rfloor}\frac{m!}{k!(m-2k)!}\left(1-\left(i\sqrt{\frac{2T}{\Delta T}}\right)^{-2}\right)^k\mathcal H_{m-2k}\left(-i\sqrt{\frac{T^*}{T^*+T}}\right)\\
	&=\left(\frac{\Delta T}{T^*+T}\right)^{m/2}\frac{i^{-m}}{m!}\mathcal H_m\left(\sqrt{\frac{2TT^*}{(T^*)^2-T^2}}\right),\quad m\in\N_0.
	\end{align*}
	Continuing \eqref{int5}, we have
	\begin{align} 
	\label{int6}
	\int_{-\infty}^\infty e^{-x_2^2/(4T)}&e^{-x_2/\sqrt{T}}\psi_m^{T^*}(x_2)\,dx_2
	\nonumber\\
	&=(8\pi T^*)^{1/4}e^{T^*/(T^*+T)}\sqrt{\frac{T}{T^*+T}}\frac{(-1)^m}{\sqrt{2^mm!}}\left(\frac{\Delta T}{T^*+T}\right)^{m/2}
	\notag\\
	&\quad\times
	\mathcal H_m\left(\sqrt{\frac{2TT^*}{(T^*)^2-T^2}}\right),\quad m\in\N_0.
	\end{align}
	Thus, using \eqref{int3} and \eqref{int6} and continuing \eqref{wnmcoef00}, we get
	\begin{align} 
	\label{wnmcoef}
	W_{nm}&=\left(\frac{\sqrt{2\pi T^*}}{e(T^*+T)}\right)^{1/2}\frac{e^{T^*/(T^*+T)}}{2^{2n+1}\sqrt{(2n+1)!}}\frac{(-1)^{m+1}}{\sqrt{2^mm!}}
	\notag\\
	&\quad\times
	\left(\frac{\Delta T}{T^*+T}\right)^{m/2}\mathcal H_m\left(\sqrt{\frac{2TT^*}{(T^*)^2-T^2}}\right), \quad
	n\in\N_0,\ m\in\N_0.
	\end{align}
	Taking into account \eqref{gpm}, \eqref{hpn}, and \eqref{wnmcoef}, we obtain
	\begin{align}
	\label{gpmexamp}
	g_{pm}^N&=\left(\frac{\sqrt{2T^* T^3}}{e(T^*+T)}\right)^{1/2}e^{T^*/(T^*+T)}\frac{(-1)^m}{\sqrt{2^mm!}}\left(\frac{\Delta T}{T^*+T}\right)^{m/2}
	\notag\\
	&\quad
	\times\mathcal H_m\left(\sqrt{\frac{2TT^*}{(T^*)^2-T^2}}\right)\frac{(-1)^p(8T)^p}{(2p+1)!}\sum_{n=p}^N \frac{1}{2^{3n}(n-p)!},
		\notag\\
	&\kern40ex   
	p=\overline{0,N},\ m=\overline{0,M}.
	\end{align} 
	
	For  $l\in\mathbb{N}$, let us construct  a control $u_{NM}^l$ in the form \eqref{umnl} with $g_{pm}^N$ defined by \eqref{gpmexamp}.	
	Then we obtain the end state of the solution $W$ to \eqref{eq1}, \eqref{ic1} with $u=u_{NM}^l$:
	\begin{equation}
	\label{apprendst}
	W(\cdot,T)=	\mathcal W_0(\cdot,T)+\mathcal W_{u_{NM}^{l}}(\cdot,T),
	\end{equation}
	where $\mathcal W_0(\cdot,T)$ is given by \eqref{w0ex}, $\mathcal W_{u_{NM}^{l}}$ is determined by \eqref{wu}.       
	
	Let us obtain estimate \eqref{lll9} in explicit form. 
	We have
	\begin{equation}
	\label{vt0}
	V^T_0=\sum_{n=0}^\infty \sum_{m=0}^\infty W_{nm}\widehat\Theta_{nm}^{TT^*}=\sum_{m=0}^\infty \widehat\omega_m\big((\cdot)_{[1]}\big)\widehat\psi_m^{T^*}\big((\cdot)_{[2]}\big),
	\end{equation}
	where
	\begin{equation}
	\label{omega}
	\widehat\omega_m\big((\cdot)_{[1]}\big)=\sum_{n=0}^\infty W_{nm}\widehat\psi_{2n+1}^T\big((\cdot)_{[1]}\big),\quad m\in\N_0.
	\end{equation}
	Substituting \eqref{hermite}, \eqref{psiext}, and \eqref{wnmcoef} into \eqref{omega}, we obtain
	\begin{align*}
	\widehat\omega_m(\sigma_1)&=i\left(\frac{2\sqrt{TT^*}}{T^*+T}\right)^{1/2}e^{T^*/(T^*+T)-3/8}\frac{(-1)^m}{\sqrt{2^mm!}}\left(\frac{\Delta T}{T^*+T}\right)^{m/2}
	\\
	&\quad
	\times\mathcal H_m\left(\sqrt{\frac{2TT^*}{(T^*)^2-T^2}}\right)e^{-T\sigma_1^2}\sin\left(\sqrt{T}\sigma_1\right),\quad m\in\N_0,\ \sigma_1\in\mathbb{R}.
	\end{align*} 
	Since
	\begin{equation*}
	\int_{-\infty}^\infty \left|e^{-T\sigma_1^2}\sin\left(\sqrt{T}\sigma_1\right)\right|^2\,d\sigma_1 \leq \sqrt{\frac{\pi}{2T}},
	\end{equation*}
	we have
	\begin{align}
	\label{normomega1}
	\left\|\widehat\omega_m\right\|_0
	&\leq\left(\frac{\sqrt{2\pi T^*}}{T^*+T}\right)^{1/2}e^{T^*/(T^*+T)-3/8}\frac{1}{\sqrt{2^mm!}}\left(\frac{\Delta T}{T^*+T}\right)^{m/2}
	\notag\\
	&\quad \times
	\mathcal H_m\left(\sqrt{\frac{2TT^*}{(T^*)^2-T^2}}\right),\quad m\in\N_0.
	\end{align}
	Let $V_{NM}$ be determined by \eqref{vtnm}. Taking into account \eqref{vt0}, we get
	\begin{align}
	\label{estim1}
	\left\|V^T_0-V_{NM}\right\|_0
	&\leq \left\|\sum_{m=M+1}^\infty \widehat\omega_m\big((\cdot)_{[1]}\big)\widehat\psi_m^{T^*}\big((\cdot)_{[2]}\big)\right\|_0
	\notag\\
	&\quad
	+\left\|\sum_{m=0}^M \left( \widehat\omega_m\big((\cdot)_{[1]}\big)-\widehat\omega_m^N\big((\cdot)_{[1]}\big)\right)\widehat\psi_m^{T^*}\big((\cdot)_{[2]}\big)\right\|_0,
	\end{align}
	where $\widehat\omega_m^N$ is determined by \eqref{omegan}.
	Taking into account the following estimate (\cite[18.14.9]{NIST})
	\begin{equation*}
	\frac{1}{\sqrt{2^nn!}}e^{-x^2/2}\left|\mathcal H_n(x)\right|\leq 1,\quad x\in\mathbb{R},\ n\in\N_0,
	\end{equation*}
	we conclude that
	\begin{equation}
	\label{estimhhh}
	\left|\mathcal H_m\left(\sqrt{\frac{2TT^*}{(T^*)^2-T^2}}\right)\right|\leq\sqrt{2^mm!}e^{TT^*/\left((T^*)^2-T^2\right)},\quad m\in\N_0.
	\end{equation}
	Using \eqref{normomega1} and \eqref{estimhhh} we obtain
	\begin{align}
	\label{estim2}
	&\left\|\sum_{m=M+1}^\infty \widehat\omega_m\big((\cdot)_{[1]}\big)\widehat\psi_m^{T^*}\big((\cdot)_{[2]}\big)\right\|_0
	\notag\\
	&\quad
	\leq \left(\frac{\sqrt{2\pi T^*}}{T^*+T}\right)^{1/2}e^{T^*/(T^*+T)-3/8}
	\notag\\
	&\qquad\times
	\left(\sum_{m=M+1}^\infty \frac{1}{2^mm!}\left(\frac{\Delta T}{T^*+T}\right)^m\left(\mathcal H_m\left(\sqrt{\frac{2TT^*}{(T^*)^2-T^2}}\right)\right)^2\right)^{1/2}
	\nonumber\\
	&\quad
	\leq \left(\frac{\sqrt{2\pi T^*}}{T^*+T}\right)^{1/2}e^{(T^*)^2/\left((T^*)^2-T^2\right)-3/8}\left(\sum_{m=M+1}^\infty\left(\frac{\Delta T}{T^*+T}\right)^m\right)^{1/2}\nonumber\\
	&\quad
	=\left(\frac{\sqrt{2\pi T^*}}{T^*+T}\right)^{1/2}e^{(T^*)^2/\left((T^*)^2-T^2\right)-3/8}
	\notag\\
	&\qquad\times
	\left(\frac{\Delta T}{T^*+T}\right)^{(M+1)/2}\left(\frac{1}{1-\Delta T/(T^*+T)}\right)^{1/2}\nonumber\\
	&\quad
	=\left(\frac{\pi T^*}{2}\right)^{1/4}\frac{e^{(T^*)^2/\left((T^*)^2-T^2\right)-3/8}}{\sqrt{T}}\left(\frac{\Delta T}{T^*+T}\right)^{(M+1)/2}.
	\end{align}
	We have taken into account that $T^*>T$, hence $\Delta T/(T^*+T)<1$.	
	
	With regard to \eqref{omegan}, \eqref{wnmcoef}, \eqref{omega}, and \eqref{estimhhh}, we have
	\begin{align} 
	\label{estim6}
	&\left\|\sum_{m=0}^M \left( \widehat\omega_m\big((\cdot)_{[1]}\big)-\widehat\omega_m^N\big((\cdot)_{[1]}\big)\right)\widehat\psi_m^T\big((\cdot)_{[2]}\big)\right\|_0\leq \left(\sum_{m=0}^M\sum_{n=N+1}^\infty \left|W_{nm}\right|^2\right)^{1/2}\nonumber\\
	&\quad
	=\left(\frac{\sqrt{2\pi T^*}}{e(T^*+T)}\right)^{1/2}e^{T^*/(T^*+T)}
		\left(\sum_{m=0}^M \frac{1}{2^mm!}\left(\frac{\Delta T}{T^*+T}\right)^m
		\right.
		\notag\\
		&\qquad\times\left.
		\left(\mathcal H_m\left(\sqrt{\frac{2TT^*}{(T^*)^2-T^2}}\right)\right)^2
	\sum_{n=N+1}^\infty \frac{1}{4^{2n+1}(2n+1)!}\right)^{1/2}
	\nonumber\\
	&\quad
	\leq \left(\frac{\sqrt{2\pi T^*}}{e(T^*+T)}\right)^{1/2}e^{(T^*)^2/\left((T^*)^2-T^2\right)}\left(\sum_{m=0}^\infty \left(\frac{\Delta T}{T^*+T}\right)^m\right)^{1/2}
	\notag\\
	&\qquad\times
	\left(\sinh(1/4)-\sum_{n=0}^N \frac{1}{4^{2n+1}(2n+1)!}\right)^{1/2}
	\nonumber\\
	&\quad
	\leq \left(\frac{\sqrt{2\pi T^*}}{e(T^*+T)}\right)^{1/2}e^{(T^*)^2/\left((T^*)^2-T^2\right)}
	\notag\\
	&\qquad\times
	\left(\frac{1}{1-\Delta T/(T^*+T)}\right)^{1/2}\left(\frac{\cosh(1/4)}{4^{2N+3}(2N+3)!}\right)^{1/2}
	\nonumber\\
	&\quad
	=\left(\frac{\pi T^*}{2}\right)^{1/4}e^{(T^*)^2/\left((T^*)^2-T^2\right)}\sqrt{\frac{\cosh(1/4)}{eT}}\frac{1}{2^{2N+3}\sqrt{(2N+3)!}}.
	\end{align} 
	Then, continuing \eqref{estim1} and taking into account \eqref{estim2} and \eqref{estim6}, we get
	\begin{align}
	\label{estim7}
	&\left\|V^T_0-V_{NM}\right\|_0
	\leq \left(\frac{\pi T^*}{2}\right)^{1/4}\frac{e^{(T^*)^2/\left((T^*)^2-T^2\right)}}{\sqrt{T}}
	\notag\\
	&\qquad\times
	\left(e^{-3/8}\left(\frac{\Delta T}{T^*+T}\right)^{(M+1)/2}+\sqrt{\frac{\cosh(1/4)}{e}}\frac{1}{2^{2N+3}\sqrt{(2N+3)!}}\right).
	\end{align} 
	
	Let $V_{NM}^l$ be determined by \eqref{vtnmlser}. Then \eqref{wnml} holds. Taking into account \eqref{vtnmser} and \eqref{vtnmlser}, we have
	\begin{equation}
	\label{estim8}
	\left\|V_{NM}-V_{NM}^l\right\|_0\leq \sum_{m=0}^M\sum_{p=0}^N \left|g_{pm}^N\right|\left\|\varphi_{2p+1}-\varphi_{2p+1}^l\right\|_0.
	\end{equation}
	Let us estimate $\left\|\varphi_{2p+1}-\varphi_{2p+1}^l\right\|_0$ for $p\in\N_0$. Applying the Tailor formula for the function $((e^x-1)/x)^{p+1}$ at $x=0$, we get
	\begin{equation}
	\label{exp1}
	\left(\frac{e^x-1}{x}\right)^{p+1}= 1 + (p+1)  x  \left(\frac{e^y-1}{y}\right)^p \frac{ye^y-e^y +1}{y^2},\quad x\geq 0,
	\end{equation}
	for some $y\in(0,x)$. Using again the Tailor formula for the function $ye^y-e^y$ at $y=0$, we obtain
	$$
	ye^y-e^y=-1+\frac{y^2}2(1+z)e^z, \quad y\geq0,
	$$
	for some $z\in(0,y)$. Since $|(1+z)e^{-z}|\leq1$, $z\geq0$, we have
	$$
	\left|\frac{ye^y-e^y +1}{y^2}\right|\leq\frac12 e^{2y},\quad y\geq0.
	$$
	Since $(e^y-1)/y\leq e^y\leq e^x$, $x\geq0$, continuing \eqref{exp1}, we get
	$$
	\left|\left(\frac{e^x-1}{x}\right)^{p+1}-1\right|
	\leq \frac12 (p+1) x e^{(p+2)x},\quad x\geq0.
	$$
	Then, for $l\geq 2(p+2)/T$, we have
	\begin{align}
	\label{estim9-a}
	\left(\left\|\varphi_{2p+1}-\varphi_{2p+1}^l\right\|_0\right)^2
	&\leq 2\left(\frac{p+1}{2l}\right)^2  
	\int_0^\infty \left(\lambda^2\right)^{2p+3} e^{-2\left(T-2(p+1)/l\right)\lambda^2}d\lambda
	\notag\\
	&\leq \left(\frac{p+1}{2l}\right)^2 2 
	\int_0^\infty \left(\lambda^2\right)^{2p+3} e^{-T\lambda^2}d\lambda
	\notag\\ 
	&=-2\left(\frac{p+1}{2l}\right)^2  \left(\frac{d}{dT}\right)^{2p+3}
	\int_0^\infty  e^{-T\lambda^2}d\lambda
	\notag\\
	&=\sqrt \pi \left(\frac{p+1}{2l}\right)^2 \frac{(4p+5)!!}{2^{2p+3}T^{2p+7/2}}.
	\end{align}
	Let us estimate $\sqrt{(4p+5)!!}$. Taking into account the Stirling formula (see, e.g. \cite[6.1.38]{HB}):
	\begin{equation}
	\label{stir}
	\sqrt{2\pi}n^{n+1/2} e^{-n} \leq n! \leq e n^{n+1/2} e^{-n}, \quad n\in\N_0,
	\end{equation}
	we obtain
	\begin{align*}
	(4p+5)!!&=\frac{(4p+5)!!(4p+6)!!}{(4p+6)!!}
	=\frac{(4p+6)!}{2^{2p+3}(2p+3)!}
	\\
	&\leq \frac{e2^{-(2p+3)}}{\sqrt{\pi}} \frac{(4p+6)^{4p+6}e^{-(4p+6)}}{(2p+3)^{2p+3}e^{-(2p+3)}}
	=\frac e{\sqrt{\pi}} 2^{2p+3} (2p+3)^{2p+3} e^{-(2p+3)}.
	\end{align*}
	Therefore,
	\begin{align*}
	\sqrt{(4p+5)!!}& \leq \frac {e^{1/2}}{\pi^{1/4}} 2^{p+3/2}(2p+3)^{p+3/2} e^{-(p+3/2)}
	\\
	&= \frac{2^{2p+3}}{\pi^{1/4}}(p+1)^{p+1+1/2} e^{-(p+1)}\left(1+\frac1{2p+2}\right)^{p+3/2}.
	\end{align*}
	Applying again the Stirling formula (see \eqref{stir}), we get
	\begin{equation}
	\label{ff}
	\sqrt{(4p+5)!!} \leq \left(\frac e\pi\right)^{3/4} 2^{2p+5/2} (p+1)!
	\end{equation}
	because
	\begin{align*}
	\left(1+\frac1{2p+2}\right)^{p+3/2}
	&=\exp\left(\left(p+\frac32\right)\ln\left(1+\frac1{2p+2}\right)\right)
	\\
	&\leq \exp\left( \frac{p+3/2}{2(p+1)}\right)\leq e^{3/4}.
	\end{align*}

	\begin{figure}[!h] 
		\begin{center}
			\begin{subfigure}[t]{0.49\linewidth}
				\centering \includegraphics[height=40mm]{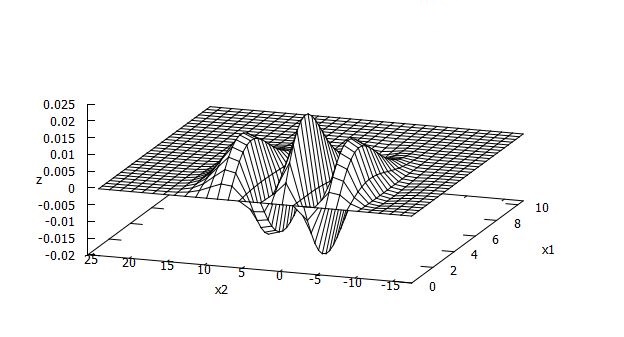}\\
				\centering \caption{\parbox[t]{0.88\textwidth}{ $N=M=3$, $l=10$.}}
				\label{fig:appr11}
			\end{subfigure}
			\begin{subfigure}[t]{0.49\linewidth}
				\centering \includegraphics[height=40mm]{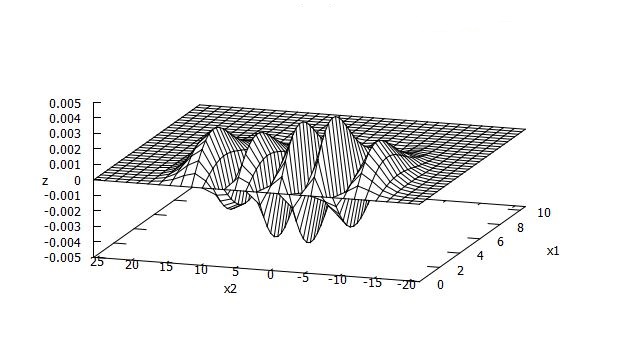}
				\centering \caption{\parbox[t]{0.88\textwidth}{%
						$N=M=6$, $l=200$.}}
				\label{fig:appr12}
			\end{subfigure}
			\centering \caption{The difference $W^T- \left(\mathcal W_0(\cdot,T)+\mathcal W_{u_{NM}^{l}}(\cdot,T)\right)$ with $T=2$ and $T^*=6$.}
			\label{fig:appr1}
		\end{center}
	\end{figure}
	\begin{figure}[!h]
		\begin{center}
			\begin{subfigure}[t]{0.49\linewidth}
				\centering \includegraphics[height=50mm]{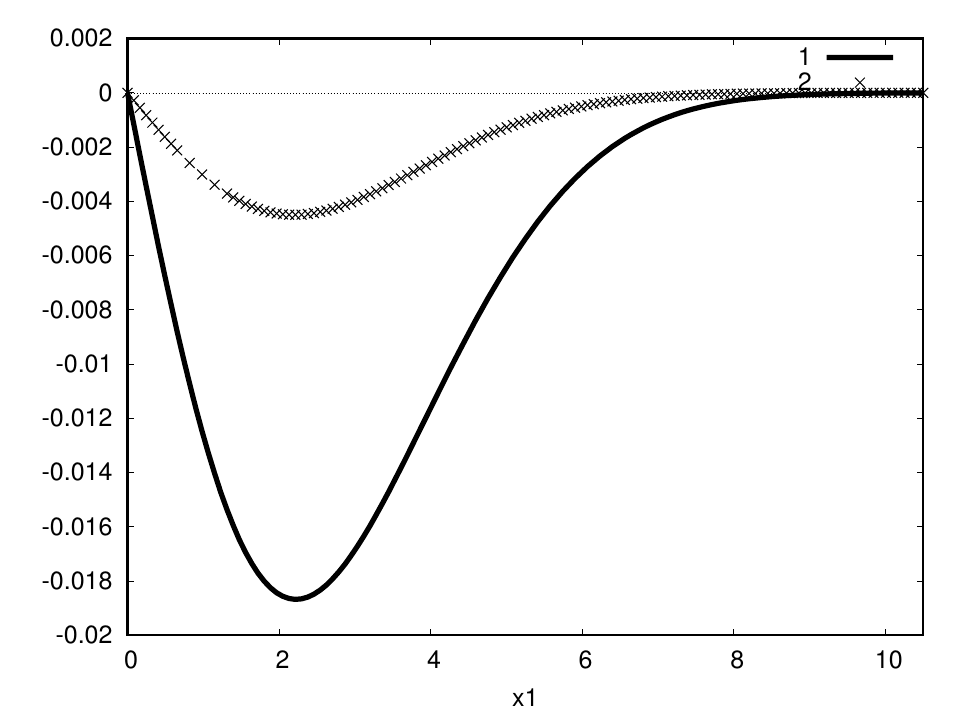}\\
				\centering \caption{\parbox[t]{0.87\textwidth}{%
						1) $N=M=3$, $l=10$;  
						2) $N=M=6$, $l=200$.}}
				\label{fig:appr21}
			\end{subfigure}
			\begin{subfigure}[t]{0.49\linewidth}
				\centering \includegraphics[height=50mm]{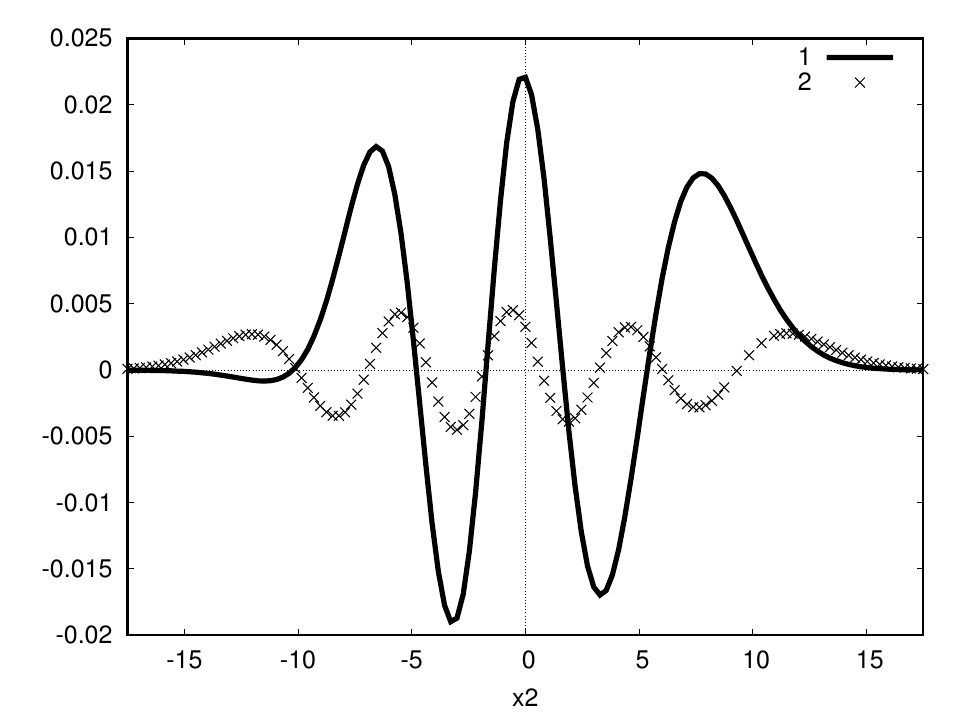}
				\centering \caption{\parbox[t]{0.87\textwidth}{%
						1) $N=M=3$, $l=10$; 
						2) $N=M=6$, $l=200$.}}
				\label{fig:appr22}
			\end{subfigure}
			\centering \caption{The difference $W^T- \left(\mathcal W_0(\cdot,T)+\mathcal W_{u_{NM}^{l}}(\cdot,T)\right)$  with $T=2$ and $T^*=6$ (vertical section for $x_2=-3$ and horizontal section for $x_1=2.2$).}
			\label{fig:appr2}
		\end{center}
	\end{figure}

	\begin{figure}[!h]
		\begin{center}
			\begin{subfigure}[t]{0.49\linewidth}
				\centering \includegraphics[height=50mm]{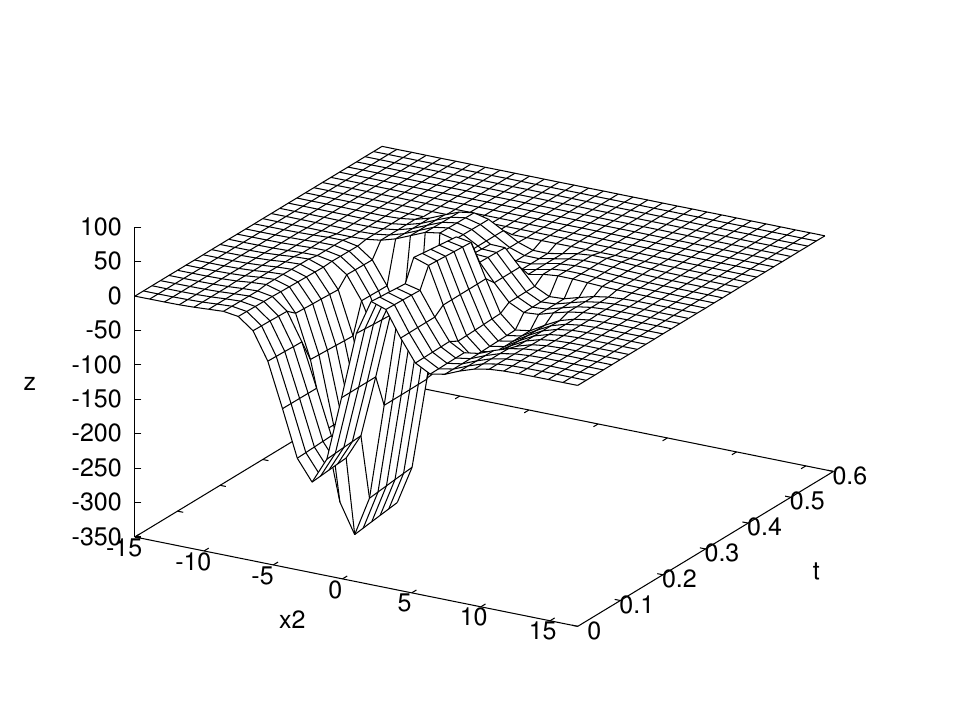}\\
				\centering \caption{\parbox[t]{0.88\textwidth}{$N=M=3$, $l=10$}}
				\label{fig:contr11}
			\end{subfigure}
			\begin{subfigure}[t]{0.49\linewidth}
				\centering \includegraphics[height=50mm]{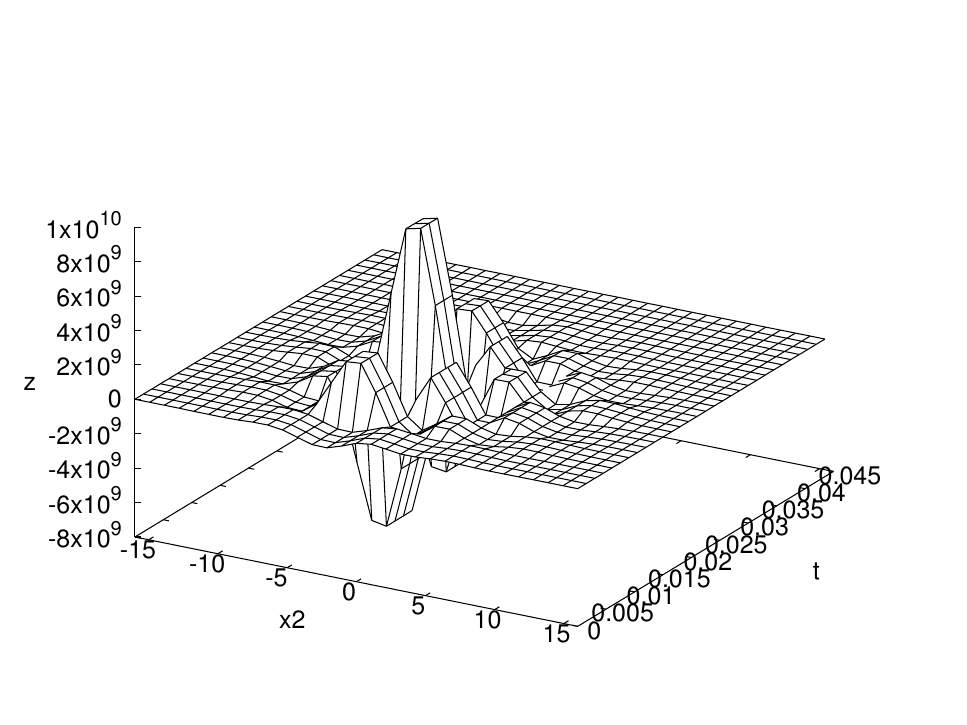}
				\centering \caption{\parbox[t]{0.88\textwidth}{%
						$N=M=6$, $l=200$}}
				\label{fig:contr22}
			\end{subfigure}
			\centering \caption{The controls $u_{NM}^{l}$ defined by \eqref{umnl} and \eqref{gpmexamp}  with $T=2$ and $T^*=6$.}
			\label{fig:contr}
		\end{center}
	\end{figure}
	
	\noindent%
	Continuing \eqref{estim9-a} and taking into account \eqref{ff}, we obtain 
	\begin{equation}
	\label{estim9}
	\left\|\varphi_{2p+1}-\varphi_{2p+1}^l\right\|_0
	\leq\frac{e^{3/4}}{\pi^{1/2}} \frac{2^p (p+1) (p+1)!}{l T^{p+7/4}},\quad p\in\N_0.
	\end{equation}
	Using \eqref{gpmexamp} and \eqref{estimhhh}, we get
	\begin{align}
	\left|g_{pm}^N\right|&\leq \left(\frac{\sqrt{2 T^* T^3}}{e(T^*+T)}\right)^{1/2}e^{T^*/(T^*+T)}\left(\frac{\Delta T}{T^*+T}\right)^{m/2}
	\notag\\
	&\quad\times
	e^{TT^*/\left((T^*)^2-T^2\right)}\frac{(8T)^p}{(2p+1)!}\sum_{j=0}^\infty \frac{1}{8^{p+j}j!}
	\nonumber\\
	&=\left(\frac{\sqrt{2 T^* T^3}}{T^*+T}\right)^{1/2}e^{(T^*)^2/\left((T^*)^2-T^2\right)-3/8}
	\notag\\
	&\quad\times
	\left(\frac{\Delta T}{T^*+T}\right)^{m/2}\frac{T^p}{(2p+1)!},\quad  p=\overline{0,N},\ m=\overline{0,M}.\label{estim10}
	\end{align}
	Therefore for $l\geq {2(N+2)}/T$, continuing \eqref{estim8} and taking into account \eqref{estim9} and \eqref{estim10}, we obtain
	\begin{align}
	\label{estim11}
	\left\|V_{NM}- V_{NM}^l\right\|_0
	&\leq \sqrt{\frac{2}{\pi(T^*+T)}}\frac{(T^*)^{1/4}e^{3/8}}{Tl}e^{(T^*)^2/\left((T^*)^2-T^2\right)}
	\notag\\
	&\quad\times
	\sum_{m=0}^M\left(\frac{\Delta T}{T^*+T}\right)^{m/2}\sum_{p=0}^N\frac{2^p(p+1)(p+1)!}{(2p+1)!}.
	\end{align}
	For the first sum in \eqref{estim11}, we have	
	\begin{align}
	\label{estim5}
	\sum_{m=0}^M\left(\frac{\Delta T}{T^*+T}\right)^{m/2}
	\leq \sum_{m=0}^\infty\left(\sqrt{\frac{\Delta T}{T^*+T}}\right)^m=\frac{\sqrt{T^*+T}}{\sqrt{T^*+T}-\sqrt{\Delta T}}.
	\end{align}
	For the second sum in \eqref{estim11}, taking into account that $(2p+1)! \geq \Big((2p)!!\Big)^2=\Big(2^p p!\Big)^2$, we have
	\begin{align}
	\label{estim111}
	\sum_{p=0}^N&\frac{2^p(p+1)(p+1)!}{(2p+1)!}
	\leq\sum_{p=0}^\infty \frac{(p+1)^2}{2^pp!}
	=\sum_{p=0}^\infty\frac{(p-1)p+3p+ 1}{2^p p!}
	\notag\\ 
	&
	=  
	\frac14 \sum_{p=2}^\infty \frac1{2^{p-2}(p-2)!} 
	+ \frac32\sum_{p=1}^\infty \frac1{2^{p-1}(p-1)!} 
	+\sum_{p=0}^\infty \frac1{2^p p!}
	=\frac{11}4 e^{1/2}.
	\end{align}
	Substituting \eqref{estim5} and \eqref{estim111} in  \eqref{estim11}, we get
	\begin{equation}
	\label{estim12}
	\left\|V_{NM}-V_{NM}^l\right\|_0
	\leq \frac{1}{l} \frac{11(T^*)^{1/4}e^{7/8}}{2T\sqrt{2\pi}}\frac{e^{(T^*)^2/\left((T^*)^2-T^2\right)}}{\sqrt{T^*+T}-\sqrt{\Delta T}},\quad l\geq \frac{2(N+2)}T.
	\end{equation}
	Taking into account \eqref{estim7}, \eqref{wnml}, and \eqref{estim12} for $l\geq 2(N+2)/T$, we have 
	\begin{align*}
	&\left\|W^T-\left(\mathcal W_0(\cdot,T)+\mathcal W_{u_{NM}^{l}}(\cdot,T)\right)\right\|^0
	\leq \left\|W_0^T-\mathcal W_{u_{NM}^{l}}(\cdot,T)\right\|^0
	\\
	&\quad
	\leq \left\|V^T_0-V_{NM}\right\|_0+\left\|V_{NM}-V_{NM}^l\right\|_0\\
	&\quad
	\leq \left(\frac{\pi T^*}{2}\right)^{1/4}\frac{e^{(T^*)^2/\left((T^*)^2-T^2\right)}}{\sqrt{T}}\left(e^{-3/8}\left(\frac{\Delta T}{T^*+T}\right)^{(M+1)/2}\right.\\
	&\qquad
	+\left.
	\sqrt{\frac{\cosh(1/4)}{e}}\frac{1}{2^{2N+3}\sqrt{(2N+3)!}}
	+\frac{1}{l}\frac{11 e^{7/8}}{\sqrt{T}(2^5\pi^3)^{1/4}\left(\sqrt{T^*+T}-\sqrt{\Delta T}\right)}\right).
	\end{align*}
	
	Thus, we have constructed the control $u_{NM}^{l}$, determined by \eqref{umnl} and \eqref{gpmexamp}, and the approximate end state, determined by \eqref{apprendst}. Let us recall that $T=2$. Put $T^*=6$. Then, we have 
	\begin{align*}
	&\left\|W^T-\left(\mathcal W_0(\cdot,T)+\mathcal W_{u_{NM}^{l}}(\cdot,T)\right)\right\|^0
	\\
	&\quad
	\leq
	\frac{(3\pi)^{1/4}}{\sqrt{2}}\left(e^{3/4}\left(\frac{1}{2}\right)^{M/2}+\frac{\sqrt{\cosh(1/4)}e^{5/8}}{2^{2N+3}\sqrt{(2N+3)!}}+\frac{1}{l}\frac{11 e^2}{(2^{11}\pi^3)^{1/4}\left(\sqrt{2}-1\right)}\right)\\
	&\quad
	=2.622851155438146\left(\frac{1}{2}\right)^{M/2}+\frac{2.350732202502537}{2^{2N+3}\sqrt{(2N+3)!}}+\frac{15.31493739172921}{l},
	\\
	&\kern60ex
	l\geq\frac{2(N+2)}{T}.
	\end{align*}
	The difference $W^T-\left(\mathcal W_0(\cdot,T)+\mathcal W_{u_{NM}^{l}}(\cdot,T)\right)$  is shown in Figure \ref{fig:appr1} and Figure \ref{fig:appr2} for the cases of $N=M=3$, $l=10$  and $N=M=6$, $l=200$. The controls $u_{NM}^{l}$ are given in Figure \ref{fig:contr} for the same cases.
	
\end{example}



\subsection*{Acknowledgments.} 
The authors are partially supported  by the Akhiezer Foundation and ``Pauli Ukraine Project'', funded in the WPI Thematic Program ``Mathematics--Magnetism--Materials'' (2024/2025).




\EndPaper



\begin{thebibliography}{99}

	\bibitem{A} 
J. Apraiz, \emph{Observability Inequalities for Parabolic Equations over Measurable Sets and Some Applications Related to the Bang-Bang Property for Control Problems,} Appl. Math. Nonlin. Sci. \textbf{2} (2017), No. 2, 543--558.

\bibitem{Barb} 
V. Barbu, \emph{Exact null internal controllability for the heat equation on unbounded convex domains,} ESAIM Control Optim. Calc. Var. \textbf{20}  (2014),  222--235.

\bibitem{CabMenZua} 
V.R. Cabanillas, S.B. De Menezes, and E. Zuazua, \emph{Null controllability in unbounded domains for the semilinear heat equation with nonlinearities involving gradient terms,} J. Optim. Theory Appl. \textbf{110} (2001), 245--264.

\bibitem{CMV} 
P. Cannarsa, P. Martinez, and  J. Vancostenoble, \emph{Null controllability of the heat equation in unbounded domains by a finite measure control region,} ESAIM Control Optim. Calc. Var. \textbf{10} (2004), 381--408.

\bibitem{LDCK}
L. Djomegne and C. Kenne, \emph{Hierarchical exact controllability of a parabolic equation with
boundary controls,} J. Math. Anal. Appl. \textbf{542} (2025), Paper No. 128799.

\bibitem{DWZh} 
Y. Duan, L.Wang, and C. Zhang, \emph{Observability inequalities for the heat equation with bounded potentials on the whole space,} SIAM J. Control Optim. \textbf{58} (2020), 1939--1960.

\bibitem{Em1} 
O.Yu. Emanuilov, \emph{Boundary controllability of parabolic equations,} Russian Math. Surveys \textbf{48} (1993), 192--194.

\bibitem{Em2} 
O.Yu. Emanuilov, \emph{Controllability of parabolic equations,} Sb. Math. \textbf{186}  (1995), 879--900.

\bibitem{FLVJMPAG05}
L.V. Fardigola, \emph{On controllability problems for the wave equation on a half-plane,} J. Math. Phys. Anal. Geom. \textbf{1} (2005,) 93--115.

\bibitem{FLVJMPAG15}
L.V. Fardigola, \emph{Modified Sobolev spaces in controllability problems  for the wave equation on a half-plane,} J. Math. Phys. Anal. Geom. \textbf{11} (2015), 18--44.

\bibitem{LVF0}
L.V. Fardigola, 
\emph{Transformation Operators and Influence Operators in Control Problems,} Thesis (Dr. Hab.), Kharkiv, 2016 (Ukrainian).

\bibitem{FKh} 
L. Fardigola and K. Khalina, \emph{Reachability and controllability problems for the heat equation on a half-axis,} J. Math. Phys. Anal. Geom. \textbf{15} (2019), 57--78.

\bibitem{FKh2} 
L. Fardigola and K. Khalina, \emph{Controllability problems for the heat equation
on a half-axis with a bounded control in the
Neumann boundary condition,} Math. Control Relat. Fields \textbf{11} (2021),  211--236.

\bibitem{FKh3} 
L. Fardigola and K. Khalina, \emph{Controllability problems for the heat equation in a half-plane controlled by the Dirichlet boundary condition with a point-wise control,} J. Math. Phys. Anal. Geom. \textbf{18} (2022),  75--104.

\bibitem{FKh4} 
L. Fardigola and K. Khalina, \emph{Controllability problems for the heat equation with variable coefficients on a half-axis,} ESAIM Control Optim. Calc. Var. \textbf{28} (2022), Paper No. 41.

\bibitem{FKh5} 
L. Fardigola and K. Khalina, \emph{Controllability problems for the heat equation with variable coefficients on a half-axis controlled by the Neumann
boundary condition,} J. Math. Phys. Anal. Geom. \textbf{19} (2023),  616--641.

\bibitem{FKh6} 
L. Fardigola and K. Khalina, \emph{Controllability problems for the heat equation in a half-plane controlled by the Neumann boundary condition with a point-wise control,} J. Math. Phys. Anal. Geom. \textbf{21} (2025),  23--55.

\bibitem{FatRus} 
H.O. Fattorini and D.L. Russell, \emph{Exact controllability theorems for linear parabolic equations in one space dimension,}  Arch. Rational Mech. Anal. \textbf{43}  (1971), 272--292.

\bibitem{FatRus1} 
H.O. Fattorini and D.L. Russell, \emph{Uniform bounds on biorthogonal functions for real exponentials with an application to the control theory of parabolic equations,}  Quart. Appl. Math. \textbf{32}   (1974/75), 45--69.

\bibitem{Fat} 
H.O. Fattorini, \emph{Boundary control of temperature distributions in a parallelepipedon,}  SIAM J. Control \textbf{13} (1975,) 1--13.

\bibitem{AFOYI}
A. Fursikov and O.Yu. Imanuvilov, \emph{Controllability of Evolution Equations,} Lecture Notes
Series 34, Research Institute of Mathematics, Global Analysis Research Center, Seoul
National University, 1996.

\bibitem{VG}
S.G. Gindikin and L.R. Volevich, \emph{Distributions and Convolution Equations,} Gordon and Breach Sci. Publ., Philadelphia, 1992.

\bibitem{WL} 
W. Gong and B. Li, \emph{Improved error estimates for semidiscrete finite element solutions of parabolic Dirichlet boundary control problems,} IMA J. Numer. Anal. \textbf{40}  (2020), 2898--2939.

\bibitem{GonTer} 
M. Gonz\'alez-Burgos and L. de Teresa, \emph{Some results on controllability for linear and nonlinear heat equations in unbounded domains,} Adv. Differential Equations \textbf{12}   (2007), 1201--1240.

\bibitem{HB}
\emph{Handbook of Mathematical Functions with Formulas Graphs and Mathematical Tables,} \emph{Eds. M. Abramowitz and I.A. Stegun,} National Bureau of Standards, Applied Mathematics Series, \textbf{55}, Washington, DC, 1972.


\bibitem{IY} 
O.Yu. Imanuvilov and M. Yamamoto, \emph{Carleman inequalities for parabolic equations
in Sobolev spaces of negative order and exact controllability for semilinear parabolic equations,} Publ. RIMS, Kyoto Univ. \textbf{39} (2003), 227--274.

\bibitem{LR} 
G. Lebeau and L. Robbiano, \emph{Contr\^ole exact de l'\'equation de la chaleur,} Comm. Partial Differential Equations \textbf{20}  (1995), 335--356.

\bibitem{LTZua} 
J. Loh\'eac, E. Tr\'elat, and E. Zuazua, \emph{Minimal controllability time for the heat equation under unilateral state or control constraints,} Math. Models Methods Appl. Sci. \textbf{27}  (2017), 1587--1644.

\bibitem{MenCab} 
S.B. De Menezes and V.R. Cabanillas, \emph{Null controllability for the semilinear heat equation in unbounded domains,}  Pesquimat \textbf{4} (2001), 35--54.

\bibitem{MRT} S. Micu, I. Roventa, and M. Tucsnak, \emph{Time optimal boundary controls for the heat equation,} J. Funct. Anal. \textbf{263} (2012), 25--49.

\bibitem{MZua1}
S. Micu and E. Zuazua, \emph{On the lack of null controllability of the heat equation on the half-line,} Trans. Amer. Math. Soc. \textbf{353} (2001), 1635--1659.

\bibitem{MZua2} 
S. Micu and E. Zuazua, \emph{On the lack of null controllability of the heat equation on the half-space,} Port. Math. (N.S.) \textbf{58} (2001), 1--24.

\bibitem{Mill} 
L. Miller, \emph{On the null-controllability of the heat equation in unbounded domains,}  Bull. Sci. Math. \textbf{129} (2005), 175--185.

\bibitem{NIST}
\emph{NIST Handbook of Mathematical Functions,} \emph{Eds. F.W.J.~Olver, D.W.~Lozier, R.F.~Boisvert, and  C.W.~Clark,} National Institute of Standards and Technology, Cambridge University Press, 2010.

\bibitem{DD} L.L.D. Njoukoue and G. Deugoue, \emph{Stackelberg control in an unbounded domain for a parabolic equation,} J. Nonl. Evol. Equ. Appl. \textbf{5} (2021), 95--118.

\bibitem{Rus} 
D.L. Russell, \emph{A unified boundary controllability theory for hyperbolic and parabolic partial differential equations,}  Stud. Appl. Math. \textbf{52} (1973), 189--211.

\bibitem{Schw}
L. Schwartz, \emph{Th\'{e}orie des distributions,} \textbf{I}, \textbf{II}, Hermann, Paris, 1950--1951.

\bibitem{SHB}
M. Shubin,
\emph{Invitation to Partial Differential Equations,}
Amer. Math. Soc., Providence, RI, 2020.

\bibitem{Ter} L. de Teresa, \emph{Approximate controllability of a semilinear heat equation in $\mathbb R^N$,}  SIAM J. Control Optim. \textbf{36} (1998), 2128--2147.

\bibitem{TerZua} L. de Teresa and E. Zuazua, \emph{Approximate controllability of a semilinear heat equation in unbounded domains,} Nonlinear Anal. \textbf{37}  (1999), 1059--1090.


\end{thebibliography}
\end{document}